\documentclass[preprint,12pt]{elsarticle}

\usepackage{hyperref}
\usepackage[colorinlistoftodos]{todonotes}
\setuptodonotes{inline, backgroundcolor=red!40!white}

\usepackage{float}

\usepackage{amssymb}
\usepackage{amsmath}
\usepackage{amsthm}

\usepackage{booktabs}
\usepackage[linesnumbered,ruled]{algorithm2e}

\newtheorem{theorem}{Theorem}

\newtheorem{proposition}{Proposition}
\newtheorem{corollary}{Corollary}
\newdefinition{remark}{Remark}
\newproof{pf}{Proof}

\newtheorem{example}{Example}

\begin{document}
	
	\begin{frontmatter}

		\title{
			Simple ordinary geometrically split abelian varieties over finite fields
		}

		\author{Novoselov S. A.} %
		
		\affiliation{organization={Immanuel Kant Baltic Federal University},%
			addressline={14 A. Nevskogo ul.}, 
			city={Kaliningrad},
			postcode={236016}, 
			country={Russia}}
		
		\begin{abstract}
			In this paper we provide a partial classification (necessary conditions) for all possible simple ordinary geometrically split abelian varieties over a finite field~$\mathbb{F}_q$. Let~$A$ be such an abelian variety, i.e.~$A$ is simple over~$\mathbb{F}_q$ and splits as~$A \sim B^s$ over~$\mathbb{F}_{q^d}$ for some simple abelian variety~$B$. Assume that~$d$ is a minimal extension degree over which~$A$ splits and~$B$ has a minimal field of definition (up to isogeny)~$\mathbb{F}_{q^e}\subseteq \mathbb{F}_{q^d}$. For dimension~$2$, Howe and Zhu showed that the only possible degrees~$d$ are~$2$, $3$, $4$, and~$6$. Moreover, for each such degree~$d$, they provided explicit formulae expressing the isogeny invariant, namely the characteristic polynomial of the Frobenius endomorphism~$\chi_{A,q}$, in terms of the coeffients of~$\chi_{B,q}$. We extend this result to ordinary abelian varieties of higher dimension. For dimensions~$3$--$10$, we provide complete lists of such varieties as components of the Weil restriction~$W_{\mathbb{F}_{q^d}/\mathbb{F}_q}(B)$ together with all possible triples~$(d, e, s)$ and restrictions on characteristic polynomial~$\chi_{A,q}$ of the form~$\chi_{A, q} \mid \chi_{V, q}$. Here, $V$ is a component of Weil restriction~$W_{\mathbb{F}_{q^d}/\mathbb{F}_q}(B)$ containing $A$ as an isogeny factor and having explicit expression of the characteristic polynomial~$\chi_{V,q}$ in terms of coefficients of~$\chi_{B,q^e}$. For higher dimensions, we propose an algorithm for enumeration of all possible cases for a given dimension of~$A$. We also provide an implementation of this algorithm.
		\end{abstract}

		\begin{keyword}
			abelian varieties \sep geometrically split \sep characteristic polynomials

		\end{keyword}
		
	\end{frontmatter}

		\section{Introduction}
		\label{sec:intro}
		
		Abelian varieties over a finite field have various applications. In particular, they can be used in cryptography for constructing discrete logarithm and isogeny based cryptosystems. Because of their efficient arithmetic, Jacobians of curves are usually chosen for this purpose. To instantiate a cryptosystem, we need to know the number of points on the selected abelian variety or, equivalently, the order of abelian variety as a group. Determining the number of points is a hard problem in the general case. For elliptic curves (abelian varieties of dimension $1$) there is an efficient algorithm due to Schoof, Elkies, and Atkin~\cite{Schoof1995}. For the Jacobians of genus~$2$ curves (dimension~$2$ abelian varieties) there is an algorithm of Gaudry-Schost~\cite{GaudrySchost2012}. However, computing the number of points in this setting already requires significant computational effort. For dimensions~$3$ and higher, such computations are already unreachable despite the fact that theoretically there is a polynomial time in~$\log q$ point-counting algorithm~\cite{AGS2019}. This motivated Dobson, Galbraith and Smith~\cite{DGS2020} to propose a variant of a verifiable delay function whose security relies on point-counting problem on the Jacobian of a genus~$3$ hyperelliptic curve (an abelian variety of dimension~$3$). Geometrically split abelian varieties decompose over an extension of the base field into a product of abelian varieties of lower dimension. Thus, we can reduce a point-counting problem from hard problem in high dimension to a simpler problem in lower dimension. By Tate's theorem~\cite{Tate1966} abelian varieties can be classified up to isogeny using the characteristic polynomial of the Frobenius endomorphism on an abelian variety. Isogenous abelian varieties have the same characteristic polynomial and the same number of points. The goal of this work is to describe all such polynomials for simple geometrically split ordinary abelian varieties. %
		
		\paragraph{Prior work} The paper~\cite{SinghMcGuireZaytsev2014} described characteristic polynomials for supersingular simple abelian varieties, i.e. ones that split into a product of supersingular elliptic curves over an extension of the base field. For geometrically split ordinary abelian surfaces such a classification can be found in~\cite{HoweZhu2002,Satoh2009,GuillevicVergnaud2012,ChouKani2014}. A method for constructing characteristic polynomials of ordinary simple geometrically split abelian varieties was described in~\cite{Jones2022}. It can be used for finding examples, but does not provide an explicit classification. %
		
		\paragraph{Our contribution} 
		In this work, we describe characteristic polynomials for geometrically split abelian varieties of dimension~$\geq 3$. Our work extends the results~\cite[\S4]{HoweZhu2002},\cite{FreemanSatoh2011,ChouKani2014} for dimension~$2$. More precisely, we explicitly describe the following data for any dimension~$g$ of a simple ordinary geometrically split abelian variety.
		\begin{enumerate}
			\item A list of all possible minimal splitting degrees~$d$. Here, minimality means that an abelian variety $A$ over~$\mathbb{F}_q$ is simple over all proper subfields of~$\mathbb{F}_{q^d}$ and it is split over~$\mathbb{F}_{q^d}$ as~$A \sim B^s$ for some~$s > 1$ and a simple abelian variety~$B$. For dimension~$2$, Howe and Zhu showed in~\cite[Theorem~6]{HoweZhu2002} that the only possible values of~$d$ are~$2$, $3$, $4$, and~$6$. We extend this result to dimension~$g \geq~3$. 
			\item For each such degree~$d$, we provide a component~$V$ of the Weil restriction~$W_{\mathbb{F}_{q^d}/\mathbb{F}_q}(B)$ such that~$A$ is an isogeny factor of~$V$. Each such component is selected in such a way that there is an explicit (symbolic) expression of the characteristic polynomial~$\chi_{V,q}$ in terms of the coefficients of~$\chi_{B,q^e}$. Here,~$e$ is a minimal up to isogeny base field extension where~$B$ is defined. Since we have~$\chi_{A,q} \mid \chi_{V,q}$, this provides a classification of $\chi_{A,q}$ as factors of~$\chi_{V,q}$.
			\item We provide an algorithm to enumerate all possible cases and its implementation.
		\end{enumerate}
		
		\paragraph{Organization of the paper} \S\ref{sec:notation} introduces the notation and definitions used throughout the paper. In~\S\ref{sec:prelims} we describe preliminaries and properties on geometrically split abelian varieties and the Weil restriction, which are used to construct the general algorithm. \S\ref{sec:algorithm} describes a general algorithm for enumeration of all possible simple ordinary geometrically split abelian varieties as components of the Weil restriction of~$B$. \S\ref{sec:dim3}, \ref{sec:dim4}, and~\S\ref{sec:prime:dim} present lists for special cases -- for dimension~$3$ and~$4$, and the case when the dimension is a prime number~$\geq 5$. \S\ref{sec:dim6-10} presents the enumeration results for dimensions~$6$, $8$, $9$, and~$10$. \S\ref{sec:examples} contains examples. In particular, for each case in dimension~$4$ there are explicit examples of ordinary simple geometrically split abelian varieties.
		
		\section{Notation}
		\label{sec:notation}
		
		\begin{itemize}
			\item $\operatorname{End}^0(A) := \operatorname{End}(A) \otimes \mathbb{Q}$.
			\item $\Phi_n(T)$ is the $n$th cyclotomic polynomial.
			\item $\phi(n)$ is Euler's totient function.
			\item $\pi_{A}$ is the $q$-power Frobenius endomorphism on an abelian variety~$A/\mathbb{F}_q$.
			\item $\chi_{A,q}$ is the characteristic polynomial of the $q$-power Frobenius endomorphism~$\pi_A$ on an abelian variety~$A/\mathbb{F}_q$.
			\item $W_{\mathbb{F}_{q^d}/\mathbb{F}_{q^e}}(A)$, $W_{d,e}(A)$ -- the Weil restriction of an abelian variety~$A$ defined over~$\mathbb{F}_{q^d}$ to the subfield~$\mathbb{F}_{q^e}$.
			\item $V_{d}(A) \simeq \ker\Phi_d(\pi_A)$ is a primitive subgroup of the Weil restriction~$W_{d,1}(A)$.
			\item $A \mid B$, an abelian variety~$A$ is an isogeny factor of an abelian variety~$B$.
		\end{itemize}
		
		\section{Preliminaries}
		\label{sec:prelims}
		\subsection{Geometrically split abelian varieties}
		\paragraph{Possible decompositions}
		In this work we study abelian varieties~$A$ over a finite field~$\mathbb{F}_q$ that are simple, but geometrically (over~$\overline{\mathbb{F}}_q$) split. For such varieties, the possible decompositions can be restricted in the following way.
		\begin{proposition}
			\label{prop:av:isotypicity}
			Let~$A$ be a simple abelian variety of dimension~$g$ over a finite field~$\mathbb{F}_q$. If $A$ is not simple over~$\mathbb{F}_{q^d}$ for some~$d>1$, then~$A \sim B^s$ over~$\mathbb{F}_{q^d}$, where~$B$ is a simple abelian variety of dimension~$\frac{g}{s}$. In particular, if $g$ is prime, then~$A \sim E^g$ over~$\mathbb{F}_{q^d}$, where~$E$ is an elliptic curve.
		\end{proposition}
		\begin{proof}
			Since $A$ is a simple abelian variety, it is trivially isotypic (i.e. a power of a simple abelian variety). It is known that for finite fields this property is preserved over all extensions of the base field~\cite[Prop.~1.2.6.1]{CCO2013}. So, if $A$ is not simple over~$\mathbb{F}_{q^d}$ for some~$d > 1$, then~$A \sim B^s$ over the field~$\mathbb{F}_{q^d}$ as claimed.
		\end{proof}
		Therefore, the cases~$A \sim B_1 \times B_2$ over~$\mathbb{F}_{q^d}$ for $B_1 \not\sim B_2$ and so on are impossible when $A$ is simple. However, this property is special for finite fields and this property need not hold over other fields.
		
		\paragraph{Possible splitting degrees}
		For each dimension, there are only finitely many base field extensions where an abelian variety can be split. We derive the degrees of such extensions from the work of Howe and Zhu~\cite{HoweZhu2002} for each given dimension~$g$. Let~$d$ be a minimal number such that a geometrically split abelian variety~$A/\mathbb{F}_q$ of dimension~$g$ splits over~$\mathbb{F}_{q^d}$. In other words,~$A$ is simple over all proper subfields of~$\mathbb{F}_{q^d}$. We call such a number~$d$ the minimal splitting degree. Now consider limitations for the degree~$d$ for ordinary abelian varieties.
		
		\begin{proposition}
			\label{prop:spliting_degrees}
			Let~$A$ be a simple ordinary abelian variety over a finite field~$\mathbb{F}_q$ such that~$A \sim B^s$ over~$\mathbb{F}_{q^d}$ for some~$s > 1$, $d>1$, and a simple abelian variety~$B$. Let $\pi$ be the Frobenius endomorphism on~$A$. Then~$\mathbb{Q}(\pi^d)$ is a proper subfield of $\mathbb{Q}(\pi)$ and there are the following two possibilities.
			\begin{enumerate}
				\item The characteristic polynomial~$\chi_{A,q}(T)$ belongs to $\mathbb{Z}[T^d]$, and
				\[
				d \in \{ z \in \mathbb{Z}_{>1} \mid z \text{~divides~} g \}.
				\]
				\item $\mathbb{Q}(\pi) = \mathbb{Q}(\pi^d, \zeta_d)$, where $\zeta_d$ is a primitive~$d$th root of unity in $\mathbb{Q}(\pi)$. Moreover,
				$d \in \{z \in \mathbb{Z}_{> 2} \mid s \leq \phi(z) \leq 2 g,  \phi(z) \textit{~divides~} 2 g, s \textit{~divides~} \phi(z) \}$.
			\end{enumerate}
		\end{proposition}
		\begin{proof}
			The abelian variety~$A$ is not simple over~$\mathbb{F}_{q^d}$, so $\mathbb{Q}(\pi^d) \neq \mathbb{Q}(\pi)$, otherwise it must be simple by Lemma 4 from~\cite{HoweZhu2002}. Therefore,  $\mathbb{Q}(\pi^d)$ is a proper subfield of~$\mathbb{Q}(\pi)$. According to Proposition~$3$ from~\cite{HoweZhu2002} there are two cases for a geometrically split abelian variety:
			\begin{enumerate}
				\item the minimal polynomial of~$\pi$ lies in~$\mathbb{Z}[T^d]$;
				\item $\mathbb{Q}(\pi) = \mathbb{Q}(\pi^d, \zeta_d)$, where~$\zeta_d$ is a primitive $d$th root of unity and~$\mathbb{Q}(\pi^d)$.
			\end{enumerate}
			Note that the minimal polynomial is equal to the characteristic polynomial for ordinary simple abelian varieties. We now consider the degrees~$d$ of extensions that can occur in each of two cases above.
			
			1. In this case~$d$ must be a divisor of~$\deg \chi_{A,q} = 2 g$ and, since~$A$ is ordinary abelian variety, the coefficient~$a_g$ of the polynomial $\chi_{A,q}$ at $T^g$ must not be divisible by $p = \operatorname{char}\mathbb{F}_q$ (otherwise, the $p$-rank of $A$ is not equal to~$g$). Thus, we have to take~$d$ as a divisor of~$g$.
			
			2. When~$d = 2$, we have~$\deg \zeta_d = 1$ and $\mathbb{Q}(\pi) = \mathbb{Q}(\pi^d)$. This would mean that~$A/\mathbb{F}_{q^d}$ is simple, which contradicts to the conditions of our proposition. So, we take~$d > 2$. The field~$\mathbb{Q}(\zeta_d)$ is a subfield of~$\mathbb{Q}(\pi)$. Therefore, $\phi(d)$ is a divisor of~$2 g$. Since the field~$\mathbb{Q}(\zeta_d)$ is Galois and~$\mathbb{Q}(\pi) = \mathbb{Q}(\pi^d, \zeta_d)$, it holds that
			\[
			\deg\pi = \frac{\deg\pi^d \cdot \deg\zeta_d}{[\mathbb{Q}(\pi^d) \cap \mathbb{Q}(\zeta_d): \mathbb{Q}]}.
			\]
			We note that it is enough for one of the fields to be Galois for equality to hold (see Corollary 3.20 from~\cite{Milne2022}). So, the field $\mathbb{Q}(\pi^d)$ is not required to be Galois. Thus, $[\mathbb{Q}(\pi^d) \cap \mathbb{Q}(\zeta_d): \mathbb{Q}] = \phi(d) / s$. 
			Since $\mathbb{Q}$-dimension of the intersection must be integer, $s$ is a divisor of~$\phi(d)$ and~$s \leq \phi(d)$.
		\end{proof}
		
		If we assume that~$d$ is a minimal splitting degree, i.e.~$A$ is simple over all proper subfields of~$\mathbb{F}_{q^d}$, then we can obtain more restrictions on possible pairs~$(d,s)$.
		
		\begin{proposition}
			\label{prop:splitting_degrees:minimal}
			Under the assumptions of Proposition~\ref{prop:spliting_degrees}, let~$A$ be simple over all proper subfields of~$\mathbb{F}_{q^d}$. Then the following holds.
			\begin{enumerate}
				\item\label{prop:splitting_degrees:minimal:case1}
				If $\chi_{A,{q}}(T) \in \mathbb{Z}[T^d]$ then $\chi_{A,{q}}(T) = \chi_{B,q^d}(T^d)$, $s = d$, and~$d$ is prime.
				\item\label{prop:splitting_degrees:minimal:case_cyc}
				If~$\mathbb{Q}(\pi) = \mathbb{Q}(\zeta_d)$ then~$2 g = \phi(d)$ and~$s$ is odd.
			\end{enumerate}
		\end{proposition}
		\begin{proof}
			\ref{prop:splitting_degrees:minimal:case1}.
			From the proof of Lemma~$5$ in~\cite{HoweZhu2002}, we know that the minimal polynomial~$m_\pi(T)$ of~$\pi$ is equal to~$m_{\pi^d}(T^d) = m_{\pi_B}(T^d)$ under conditions of our proposition. Since $A$ and~$B$ are ordinary and simple, we have~$\chi_{A,q}(T) = \chi_{B,q}(T^d)$. So, $2 g = \deg\chi_{A,q}(T) = \deg\chi_{B,q}(T^d) = 2 g d / s$ and thus $s = d$. To prove that~$d$ is prime, we use properties of the Weil restriction (see~\S\ref{sec:wres:structure}). Let~$e$ be a positive integer such that~$1 \leq e < d$ and~$e \mid d$. Since~$A$ is defined over~$\mathbb{F}_{q^e}$, we have from Proposition~\ref{prop:av_is_factor_of_W} that~$A/\mathbb{F}_{q^e}$ is an isogeny factor of~$W = W_{\mathbb{F}_{q^d}/\mathbb{F}_{q^e}}(B)$. Therefore, $g = \dim A \leq \dim W = \frac{d}{e} \frac{g}{s}$. Hence, $e \leq \frac{d}{s} = 1$. Since the only remaining divisors of~$d$ are~$e=1$ and~$e = d$, we obtain that~$d$ is prime.
			
			\ref{prop:splitting_degrees:minimal:case_cyc}. To exclude the case in which~$s$ is even, we generalize the proof of Theorem~$6$ from~\cite{HoweZhu2002} for~$g=2$. Let~$K = \mathbb{Q}(\pi) = \mathbb{Q}(\zeta_d)$ and $L = \mathbb{Q}(\pi^d)$. Note that $\#\operatorname{Gal}(K/L) = [K : L] = \frac{\dim A}{\dim B} = s$. If~$s$ is even, then there is an automorphism~$\sigma$ of order~$2$ in~$\operatorname{Gal}(K/L)$ such that~$\sigma(\pi) = \zeta_d \pi$. Applying~$\sigma$ to the last equality, we obtain that~$\pi = \sigma(\zeta_d) \cdot \sigma(\pi) = \zeta_d \cdot \sigma(\zeta_d) \cdot \pi$. Therefore, $\sigma(\zeta_d) = 1/\zeta_d$ and~$\sigma$ is the complex conjugation, since it is the only automorphism with such a property. This means that~$L$ is totally real as the fixed field of~$\sigma$. This contradicts the fact that~$L = \mathbb{Q}(\pi^d) = \mathbb{Q}(\pi_B)$ is a CM-field. Thus, $s$ cannot be even.
		\end{proof}
		Note that there are explicit examples (see dimension~$3$ case in~\S\ref{sec:dim3}) of simple ordinary geometrically split abelian varieties~$A$ with odd~$s$ and~$\mathbb{Q}(\pi) = \mathbb{Q}(\zeta_d)$.
		
		\begin{remark}
			\label{remark:example:cyc_d} There are also examples~\cite[\href{https://www.lmfdb.org/Variety/Abelian/Fq/2/5/ac_ab}{2.5.ac\_ab}]{LMFDB}%
			 of the ordinary simple geometrically split abelian varieties with~$\mathbb{Q}(\pi) = \mathbb{Q}(\zeta_d)$ and even~$s$. But in this case~$d$ is not a minimal splitting degree.
		\end{remark}
		
		\subsubsection{The case of $\mathbb{Q}(\pi^d) \cap \mathbb{Q}(\zeta_d) \neq \mathbb{Q}$.}
		If the intersection of $\mathbb{Q}(\pi^d)$ and~$\mathbb{Q}(\zeta_d)$ is non-trivial, then we know (up to isomorphism) a subfield of~$\mathbb{Q}(\pi^d) = \mathbb{Q}(\pi_B)$ equal to this intersection~$L$. In some cases, this allows us to determine whether the abelian variety~$B$ is absolutely simple or not.
		
		\begin{proposition}
			Let~$B$ be an absolutely simple ordinary abelian variety over a finite field~$\mathbb{F}_q$. Then~$K = \operatorname{End}^0(B)$ is a primitive CM-field, i.e.~$K$ does not contain a proper CM-subfield. In particular, let~$B$ be  an abelian surface. Then~$B$ can be absolutely simple only if~$K$ does not contain an imaginary quadratic subfield.
		\end{proposition}
		\begin{proof}
			By results of Lubin-Serre-Tate~\cite[Theorem~8.10]{CohFre+2005} and Deligne~\cite[\S2.4]{JW2019} for an ordinary abelian variety~$B$ there exists the unique up to isomorphism abelian variety~$B_\mathbb{C}$ over~$\mathbb{C}$ such that~$B_\mathbb{C}$ reduces to~$B$ modulo~$p = \operatorname{char}(\mathbb{F}_q)$ and~$\operatorname{End}(B_\mathbb{C}) \simeq \operatorname{End}(B)$. Under this isomorphism the CM-field~$K = \mathbb{Q}(\pi_B)$ lifts to a CM-field~$K_\mathbb{C}$. For an ordinary absolutely simple abelian variety~$B$, the abelian variety~$B_\mathbb{C}$ is also simple. And if~$B_\mathbb{C}$ is simple, it has primitive CM-type~\cite[\S8.2]{Shimura1998} and the field~$K_\mathbb{C}$ is primitive~(follows from~\cite[Lemma~3.5]{Streng2010}). Therefore, the isomorphic field~$K$ is also primitive. For abelian surfaces it is known~\cite[Lemma~3.4]{Streng2010} that~$K_{\mathbb{C}}$ is primitive if and only if $K_\mathbb{C}$ does not contain an imaginary quadratic subfield.
		\end{proof}

		Information about the intersection~$L$ can be used to restrict the number of possible cases for geometrically split abelian varieties. Algorithm~\ref{alg:end:enum} allows us to enumerate all suitable intersections~$\mathbb{Q}(\pi^d) \cap \mathbb{Q}(\zeta_d) \neq \mathbb{Q}$. If this algorithm returns~$\emptyset$ for $\mathbb{Q}(\pi^d) \cap \mathbb{Q}(\zeta_d) \neq \mathbb{Q}$, then we can exclude a triple~$(g,s,d)$ as impossible. The algorithm enumerates all CM-fields of $\mathbb{Q}(\zeta_d)$ in case when~$\mathbb{Q}(\zeta_d) = \mathbb{Q}(\pi_A)$. For the case when we have non-trivial intersection and this degree-$m$ intersection is a real field, the field~$\mathbb{Q}(\pi_B)$ must contain a subfield of degree~$2 m$ due to the presence of complex conjugation. So, the algorithm additionally checks that~$2 m$ divides~$\deg\pi_B = 2 g/s$.
		
		\begin{algorithm}[H]
			\caption{($\operatorname{EnumEnd}$) List of possible subfields of the endomorphism ring of~$B$.}
			\label{alg:end:enum}
			\KwData{A triple $(g,s,d)$ such that~$A \sim B^s$ over~$\mathbb{F}_{q^d}$ and~$g = \dim A$.}
			\KwResult{A set of possible subfields of~$\operatorname{End}^0(B)$.
			}
			$m \gets \phi(d) / s$\;
			$S \gets \emptyset$\;
			\uIf{$m = 2 g / s$}{
				$S \gets \{ \text{all CM-subfields of $\mathbb{Q}(\zeta_d)$ of degree $m$} \}$
			}
			\ElseIf{$m > 1$}{
				\For{$M \in \{\text{degree~$m$ subfields of $\mathbb{Q}(\zeta_d)$} \}$}{
					\If{$M$ is imaginary or $M$ is totally real and $2 m \mid 2 g/s$}{
						$S \gets S \cup \{M\}$
					}
				}
			}
			\Return{$S$}
		\end{algorithm}
		
		\begin{proposition}
			\label{prop:enum:end:props}
			If there is an ordinary abelian variety $A$ of dimension~$g$ over a finite field~$\mathbb{F}_q$ such that
			\begin{enumerate}
				\item $\mathbb{Q}(\pi_A) = \mathbb{Q}(\pi_A^d, \zeta_d)$;
				\item $A \sim B^s$ over~$\mathbb{F}_{q^d}$ for some~$s>1$, $d>1$, and a simple abelian variety~$B$;
				\item $A$ is simple over all proper subfields of~$\mathbb{F}_{q^d}$
			\end{enumerate}
			then Algorithm~\ref{alg:end:enum} given a triple $(g,s,d)$ returns a set of subfields of~$\mathbb{Q}(\zeta_d)$ of degree~$\phi(d) / s$, such that~$\mathbb{Q}(\pi_B)$ always contains one of them (up to isomorphism). If the output of the algorithm is the empty set then either~$\phi(d) = s$, or there are no abelian varieties~$A$ and~$B$ satisfying the above conditions. 
		\end{proposition}
		\begin{proof}
			We have~$\mathbb{Q}(\pi_A^d) = \mathbb{Q}(\pi_B)$ since~$A$ is ordinary and~$A \sim B^s$ over~$\mathbb{F}_{q^d}$. The number~$m$ in the algorithm is the degree of intersection~$\mathbb{Q}(\zeta_d) \cap \mathbb{Q}(\pi_B)$. Since the field~$\mathbb{Q}(\zeta_d)$ is Galois, the degree of the compositum can be computed using Corollary 3.20 from~\cite{Milne2022}. We note that~$\mathbb{Q}(\pi_B)$ is not required to be Galois. The degree of intersection~$m = \deg \mathbb{Q}(\zeta_d) \cap \mathbb{Q}(\pi_B) = \phi(d) / s$. If~$m=1$, then $\phi(d) = s$ and there is no information about subfields of~$\mathbb{Q}(\pi_B)$ available. So, the algorithm returns the empty set.
			
			In the case when~$m = 2g / s$, we have~$\mathbb{Q}(\pi_B) \subset \mathbb{Q}(\zeta_d)$. Since~$B$ is an abelian variety, the field~$\mathbb{Q}(\pi_B)$ is a CM-field of degree~$2 g / s$. Therefore, $\mathbb{Q}(\pi_B)$ must be isomorphic to one of the CM-subfields of~$\mathbb{Q}(\zeta_d)$ of degree~$2g / s$. If the field~$\mathbb{Q}(\zeta_d)$ does not contain such a subfield, then we have a contradiction to the three conditions of our proposition.
			
			Now suppose that~$m>1$. In this case the intersection is not equal to~$\mathbb{Q}$. So, we have some information about a subfield~$M$ of~$\mathbb{Q}(\pi_B)$ of degree~$m$. This subfield is also a subfield of the field~$\mathbb{Q}(\zeta_d)$. So, we can enumerate subfields of~$\mathbb{Q}(\zeta_d)$ of degree~$m$. If~$M$ is real, then due to the existence of complex conjugation it holds that~$2 m = 2 \deg M \mid \deg\pi_B = 2 g / s$. Therefore, we add to output list only real fields that satisfy this condition. If~$M$ is imaginary there is no restrictions on the degree, so we just add it to the output set~$S$. If the field~$\mathbb{Q}(\zeta_d)$ does not contain required subfields then no such abelian variety~$B$ exists.
		\end{proof}
		
		\begin{example}
			Let~$g = \dim A = 6$, $d = 5$, $A \sim B^2$ over~$\mathbb{F}_{q^5}$, and~$\mathbb{Q}(\pi_A) = \mathbb{Q}(\pi_A^5, \zeta_5)$. We have the degree of intersection~$m = \phi(d) / s = 2$. The field~$\mathbb{Q}(\zeta_5)$ has only one subfield~$M$ of degree~$2$ isomorphic to~$\mathbb{Q}(\sqrt{5})$. However, we have $2 \deg M = 4 \nmid \deg \pi_B = 6$ and Algorithm~\ref{alg:end:enum} will return empty set. Therefore, there is no geometrically split ordinary abelian variety~$A$ with such $g, d$, and~$s$. 
		\end{example}
		
		\subsection{Structure of the Weil restriction}
		\label{sec:wres:structure}
		
		To determine the characteristic polynomials of the Frobenius endomorphism of a geometrically split abelian variety we use the decomposition of the Weil restriction~$W_{\mathbb{F}_{q^d}/\mathbb{F}_q}(B)$. For the definition and properties of the Weil restriction, we refer to~\cite[\S7]{CohFre+2005}. The following proposition is the key fact that allows us to use the decomposition of the Weil restriction.
		\begin{proposition}
			\label{prop:av_is_factor_of_W}
			Let~$A$ be a simple abelian variety over a finite field~$\mathbb{F}_q$ such that~$A \sim B^s$ over~$\mathbb{F}_{q^d}$ for some simple abelian variety~$B$, $s \geq 1$, and~$d > 1$. Then $A$ is isogenous over~$\mathbb{F}_q$ to a subvariety of the Weil restriction~$W_{\mathbb{F}_{q^d}/\mathbb{F}_q}(B)$.
			\begin{proof}
				The variety~$A$ can be embedded into~$W = W_{\mathbb{F}_{q^d}/\mathbb{F}_q}(A \otimes {\mathbb{F}_{q^d}})$ as described in~\cite[\S7 and Prop.~7.13]{CohFre+2005}, so~$A$ is isogenous to an abelian subvariety of~$W$. Since~$W \sim W_{\mathbb{F}_{q^d}/\mathbb{F}_q}(B)^s$ and~$A$ is simple, we conclude that~$A$ is a factor of~$W_{\mathbb{F}_{q^d}/\mathbb{F}_q}(B)$. A proof for the case in which~$B$ is defined over~$\mathbb{F}_q$ can also be found in~\cite[Prop.~3.1]{FreemanSatoh2011}. For dimension~$2$, this also appears in~\cite[\S4]{ChouKani2014}.
			\end{proof}
		\end{proposition}
		Thus,~$A$ is a factor of the Weil restriction~$W = W_{\mathbb{F}_{q^d}/\mathbb{F}_{q}}(B)$ and this means that~$\chi_{A,q}(T)$ divides the characteristic polynomial~$\chi_{W,q}(T)$ of the Frobenius endomorphism on the Weil restriction~$W$. From~\cite[Prop.~3.2]{FreemanSatoh2011} we know that~$\chi_{W,q}(T) = \chi_{B,q^d}(T^d)$. Thus, we can always say that~$\chi_{A,q}(T) \mid \chi_{B,q^d}(T^d)$. More detailed information can be obtained by studying the structure of the Weil restriction. In particular, if~$W$ is split, then~$A$ is one of its components. If~$W$ is simple, then~$A\sim W$, and~$W$ must have dimension~$\dim A$.
		
		\subsubsection{Primitive groups}
		\label{sec:char_poly_from_prim_groups}
		Further information about characteristic polynomials can be obtained from the decomposition of the Weil restriction into primitive groups. However, this decomposition works only if~$B$ has an isogeny twist defined over~$\mathbb{F}_q$. Let~$\widetilde{B}$ be such a twist, i.e. an abelian variety defined over~$\mathbb{F}_q$ such that~$\widetilde{B} \sim B$ over~$\mathbb{F}_{q^d}$. In this case we can~\cite[Cor.~3.7]{FreemanSatoh2011} decompose~$W$ as
		\begin{equation}
			\label{eq:res:decomp:prim:groups}
			W \sim \oplus_{e \mid d} V_e(B),
		\end{equation}
		where~$V_e(B) \simeq \operatorname{ker}(\Phi_e(\pi))$, $\dim V_e(B) = \phi(e)$, and~$\Phi_e$ is the $e$th cyclotomic polynomial. Since~$A$ is simple, it is a factor of one of the varieties~$V_e(B)$ with $\dim V_e(B) = \phi(e) \cdot \dim B$~$\geq g$ and~$\chi_{A,q}(T)$ is a divisor of~$\chi_{V_e(B),q}(T)$ by Tate's theorem. The characteristic polynomial~$\chi_{V_e(B),q}(T)$ itself can be found using the following result.
		\begin{proposition}
			\label{prop:Ve:charpoly}
			Let~$A$ be an abelian variety of dimension~$g$ over a finite field~$\mathbb{F}_{q}$, let~$\pi$ be the Frobenius endomorphism on~$A$ and~$V_e(A) = \ker\Phi_e(\pi)$. Then
			\[
			\chi_{V_e(A),q}(T) = \prod_{\substack{1 \leq j \leq e \\ (e,j)=1}} \zeta^{2 g j} \chi_{A,q}(T/\zeta^j),
			\]
			where~$\zeta$ is an $e$th root of unity.
		\end{proposition}
		\begin{proof}
			Let~$\lambda_1, \ldots, \lambda_{2g}$ be complex roots of~$\chi_{A,q}(T)$. Then by the formula from~\cite[Prop.~3.4]{FreemanSatoh2011} we have
			\begin{equation*}
				\chi_{V_e(A),q}(T) =
				\prod_{i=1}^{2 g} \lambda_i^{\phi(e)} \Phi_e(T / \lambda_i)
				=
				\prod_{i=1}^{2 g} \prod_{\substack{1 \leq j \leq e \\ (e,j)=1}} (T-\zeta^j \lambda_i)
				=
				\prod_{\substack{1 \leq j \leq e \\ (e,j)=1}} \zeta^{2 g j} \chi_{A,q}(T/\zeta^j).
			\end{equation*}
		\end{proof}
		The script that computes the characteristic polynomial~$\chi_{V_e(A),q}$ for primitive groups~$V_e(A)$ in terms of the coefficients of the characteristic polynomial~$\chi_{A,q}$ is available in our repository.\footnote{This is \texttt{av\_gsplit\_charpoly.sage} script available at~\url{https://github.com/novoselov-sa/avGsplit/}} 
		
		When~$d$ is minimal such that~$s > 1$, the abelian variety~$A$ must be a factor of some primitive group~$V_e$ for every subfield of~$\mathbb{F}_{q^d}$. This fact allows us to restrict the list of possible~$d$ by excluding values for which the following proposition does not hold.
		\begin{proposition}
			\label{prop:av:factor:of:Ve}
			Let~$A$ be a simple abelian variety of dimension~$g$ over a finite field~$\mathbb{F}_{q}$ and let~$A \sim B^s$ over~$\mathbb{F}_{q^d}$, where~$d>1$, $s>1$, and $B$ is a simple abelian variety. Assume that~$A/\mathbb{F}_{q^v}$ is simple for every integer~$v \mid d$, $1 \leq v < d$ and the abelian variety~$B$ is defined (up to isogeny) over~$\mathbb{F}_q$. Then for every such~$v$ %
			the abelian variety~$A/\mathbb{F}_{q^v}$ is an isogeny factor of~$V_{\frac{d}{v}}(B/\mathbb{F}_{q^v})$ and~$\phi(\frac{d}{v}) \geq s$.
		\end{proposition}
		\begin{proof}
			We have that~$A/\mathbb{F}_{q^v}$ is a factor of~$W_{d,v}(B) \sim \oplus_{e \mid \frac{d}{v}} V_e(B/\mathbb{F}_{q^v})$. Since the abelian variety~$A/\mathbb{F}_{q^v}$ is simple, it is a factor of~$V_e(B/\mathbb{F}_{q^v})$ of suitable dimension, i.e.~$\phi(e) g/s = \dim V_e(B/\mathbb{F}_{q^v})  \geq \dim A = g$. Thus, $\phi(e) \geq s$. Since~$e \mid \frac{d}{v}$ it holds that~$e v \mid d$. Then~$V_e(B/\mathbb{F}_{q^v})$ is a component of~$W_{e v, v}(B)$ and~$V_e(B/\mathbb{F}_{q^v}) \sim B^{\phi(e)}$ over~$\mathbb{F}_{q^{ev}}$. Therefore, $A/\mathbb{F}_{q^{e v}} \mid B^{\phi(e)}$. But~$A$ is simple over every proper subfield of~$\mathbb{F}_{q^d}$. This implies~$ev = d$ and~$e = \frac{d}{v}$.
		\end{proof}

		\subsubsection{The case when abelian variety~$B$ is primitive.} If $B/\mathbb{F}_{q^d}$ is primitive, i.e. there is no abelian variety~$\widetilde{B}$ over a proper subfield of~$\mathbb{F}_{q^d}$ such that~$B \sim \widetilde{B}$, then in most cases the Weil restriction of~$B$ is simple, and hence~$A \sim W$. The precise statement is as follows.
		\begin{proposition}[Theorem 2 of \cite{DiemNaumann2003}]
			\label{prop:wr:simplicity}
			Let~$B/\mathbb{F}_{q^d}$ be a simple primitive abelian variety and let~$\pi_B$ be the Frobenius endomorphism on~$B$. If one of the following conditions holds
			\begin{enumerate}
				\item $\operatorname{End}(B)$ is commutative and further, if~$4 \mid d$, then $\pi_B \not\in -4 \mathbb{Q}[\pi_B]^4$;
				\item $d$ is prime.
			\end{enumerate}
			Then $W_{\mathbb{F}_{q^d}/\mathbb{F}_q}(B)$ is simple.
		\end{proposition}
		For an ordinary abelian variety~$B$, the endomorphism algebra~$\operatorname{End}^0(B)$ is always commutative, so~$W_{\mathbb{F}_{q^d}/\mathbb{F}_q}(B)$ is always simple for~$4 \nmid d$ and primitive~$B/\mathbb{F}_{q^d}$. Therefore, we have the following result.
		\begin{corollary}
			\label{cor:wr:primitive:4nd:case}
			Let~$A$ be an ordinary simple abelian variety over a finite field~$\mathbb{F}_q$ and let~$A \sim B^s$ over~$\mathbb{F}_{q^d}$ for some integers~$s>1$ and~$d > 1$. If $B$ is simple and primitive over~$\mathbb{F}_{q^d}$, $4 \nmid d$ or $4 \mid d$ and~$\pi_B \not\in -4 \mathbb{Q}(\pi_B)^4$, then
			\begin{enumerate}
				\item $A \sim W_{\mathbb{F}_{q^d}/\mathbb{F}_{q}}(B)$ and~$d = s$;
				\item $\chi_{A,{q}}(T) = \chi_{B,{q^d}}(T^d)$.
			\end{enumerate}
		\end{corollary}

		Suppose that~$4 \mid d$, that~$\pi_B \in -4\mathbb{Q}[\pi_B]^4$, and that~$B$ is ordinary. In this case the characteristic polynomial of the Frobenius on~$W_{\mathbb{F}_{q^d}/\mathbb{F}_q}(B)$ is reducible and so~$W_{\mathbb{F}_{q^d}/\mathbb{F}_q}(B)$ cannot be simple. More precisely, we have the following result.
		
		\begin{proposition}
			\label{prop:wr:primitive:4dcase}
			Let~$B/\mathbb{F}_{q^d}$ be a simple primitive ordinary abelian variety, $\pi_B$ be the Frobenius endomorphism on~$B$, $4 \mid d$, and~$\pi_B \in -4 \mathbb{Q}[\pi_B]^4$. Then~$W_{\mathbb{F}_{q^d}/\mathbb{F}_q}(B) \sim W_0 \times \widetilde{W}_0$, where~$W_0$ is an abelian variety of dimension~$\frac{d \dim B}{2}$ and~$\widetilde{W}_0$ is a quadratic twist of~$W_0$. Moreover, if there is an abelian variety~$C$ isogenous to a proper abelian subvariety of~$W_0$ or~$\widetilde{W}_0$ then $\dim B < \dim C < (\frac{d}{2}-1) \dim B$.
		\end{proposition}
		\begin{proof}
			Let~$d = 4 d_0$ for an integer~$d_0$, $\pi_B = - 4 \cdot f^4$ for~$f \in \mathbb{Q}[\pi_B]$, and~$W = W_{\mathbb{F}_{q^d}/\mathbb{F}_q}(B)$. Then we have the isomorphism (Corollary~6 from \cite{DiemNaumann2003} applied to ordinary abelian varieties):
			\[
			\operatorname{End}^0(W) \simeq \mathbb{Q}[\pi_B][X] / (X^d - \pi_B).
			\]
			The polynomial~$X^d-\pi_B$ is reducible and factors as
			\[
			X^d - \pi_B = (X^{2 d_0} + 2 f X^{d_0} + 2 f^2) \cdot (X^{2 d_0} - 2 f X^{d_0} + 2 f^2).
			\]
			Therefore, 
			\[
			\operatorname{End}^0(W)
			\simeq
			\mathbb{Q}[\pi_B][X] / (X^{2 d_0} + 2 f X^{d_0} + 2 f^2)
			\times
			\mathbb{Q}[\pi_B][X] / (X^{2 d_0} - 2 f X^{d_0} + 2 f^2).
			\]
			Thus, the abelian variety~$W$ is not simple and
			\[
			\operatorname{End}^0(W)
			\simeq
			\operatorname{End}^0(W_0) \times \operatorname{End}^0(\widetilde{W}_0)
			\]
			for some abelian varieties~$W_0$ and~$\widetilde{W}_0$ with 
			\[
			\operatorname{End}^0(W_0)
			\simeq
			\mathbb{Q}[X]/\chi_{W_0,q}(X)
			\simeq
			\mathbb{Q}[\pi_B][X] / (X^{2 d_0} + 2 f X^{d_0} + 2 f^2)
			\]
			and
			\[
			\operatorname{End}^0(\widetilde{W}_0)
			\simeq
			\mathbb{Q}[X]/\chi_{\widetilde{W}_0,q}(X)
			\simeq
			\mathbb{Q}[\pi_B][X] / (X^{2 d_0} - 2 f X^{d_0} + 2 f^2).
			\]
			We have~$\chi_{W_0,q}(X) = \chi_{\widetilde{W}_0,q}(-X)$, so~$W$ and~$\widetilde{W}_0$ are quadratic twists.
			
			Let~$C/\mathbb{F}_q$ be an abelian variety isogenous to a subvariety of~$W_0$ or~$\widetilde{W}_0$. Assume that~$\dim C \leq \dim B$. Since~$C$ is a factor of~$W$, it is a factor of~$W \oplus \mathbb{F}_{q^d} \sim B^{d}$. So, $C$ is a factor of~$B^d$. Since~$B$ is simple and abelian varieties over finite fields are isotypic, the case~$\dim C < \dim B$ is impossible. In case of~$\dim C = \dim B$ we have a contradiction to the primitivity of~$B$, since~$C$ is defined over~$\mathbb{F}_q$. Thus, $\dim C > \dim B$ as claimed.
			
			The abelian variety~$C$ is an isogeny factor of~$W_0 \sim_{\mathbb{F}_{q^d}} B^{d/2}$ (the case when~$C$ is a factor of~$\widetilde{W}_0$ is the same). So, $W_0 \sim C \times W_0'$ for some abelian variety~$W_0'$ with~$\dim W_0' = \frac{d}{2}\dim B - \dim C$. Since~$W_0'$ is a subvariety of~$W_0$, we have~$\dim W_0' > \dim B$. Thus, $\dim C < (\frac{d}{2}-1) \dim B$
		\end{proof}
		\begin{remark}
			Let~$A$ be a simple abelian variety over~$\mathbb{F}_{q}$ such that~$A \sim B^s$ over~$\mathbb{F}_{q^d}$. We can use the condition~$\dim C > \dim B$ in the proposition above to exclude the case when~$A$ is a simple factor of~$W_{\mathbb{F}_{q^d}/\mathbb{F}_q}(B) \sim W_0 \times \widetilde{W}_0$ and~$W_0$ (or~$\widetilde{W}_0$) splits as~$W_0 \sim A \times C$ with~$\dim C \leq \dim B$. In other words, we always have either~$\dim W_0 = \dim A$ (i.e.~$d = 2 s$), or~$d > 2 s + 2$ when~$4 \mid d$, $\pi_B \in \mathbb{Q}[\pi_B]^4$, and~$B$ is primitive over~$\mathbb{F}_{q^d}$.
		\end{remark}
		
		\subsubsection{The case when abelian variety~$B$ is defined over an intermediate subfield.}
		In this section, we study the case in which the abelian variety~$B$ has an isogeny twist defined over a proper subfield of~$\mathbb{F}_{q^d}$.
		
		\begin{proposition}
			\label{prop:wr:decomposition:intermediate_field}
			Let~$A$ be a simple ordinary abelian variety of dimension~$g$ over a finite field~$\mathbb{F}_q$, such that~$A \sim B^s$ over~$\mathbb{F}_{q^d}$, where~$B$ is a simple abelian variety,~$d$ is minimal such that~$s > 1$, and there is an abelian variety~$\widetilde{B} \sim B$ defined over~$\mathbb{F}_{q^e} \subset \mathbb{F}_{q^d}$ for~$e > 1$. Assume that~$\widetilde{B}$ is primitive over~$\mathbb{F}_{q^e}$.
			Then the following holds.
			\begin{enumerate}
				\item\label{prop:wr:decomposition:intermediate_field:Wde}
				$W_{\mathbb{F}_{q^d}/\mathbb{F}_{q^e}}(B) \sim \oplus_{r \mid \frac{d}{e}} V_r(B)$.
				\item\label{prop:wr:decomposition:intermediate_field:Ae}
				The abelian variety~$A$ is isogenous over~$\mathbb{F}_{q^e}$ to a subvariety of~$V_{\frac{d}{e}}(B)$.
				\item\label{prop:wr:decomposition:intermediate_field:frobA}
				The characteristic polynomial~$\chi_{A,q}(T)$ divides~$\chi_{V_{\frac{d}{e}}(B),q^e}(T^e)$.
				\item\label{prop:wr:decomposition:intermediate_field:A}
				The abelian variety~$A$ is isogenous over~$\mathbb{F}_q$ to a subvariety of~$W_{\mathbb{F}_{q^e}/\mathbb{F}_q}(V_{\frac{d}{e}}(B))$.
				\item
				\label{prop:wr:decomposition:intermediate_field:Ae:isogV}
				If $A/\mathbb{F}_{q^e} \sim V_{\frac{d}{e}}(B)$ then $s = \phi(\frac{d}{e})$.
				\item\label{prop:wr:decomposition:intermediate_field:Ae:subvarV}
				If $A/\mathbb{F}_{q^e}$ is isogenous to a proper subvariety of~$V_{\frac{d}{e}}(B)$ then $s < \phi(\frac{d}{e})$.
				\item
				\label{prop:wr:decomposition:intermediate_field:factor_of_Vte}
				For every finite field~$\mathbb{F}_{q^v}$ satisfying~$\mathbb{F}_{q^e} \subseteq \mathbb{F}_{q^v} \subset \mathbb{F}_{q^d}$, the abelian variety~$A/\mathbb{F}_{q^v}$ is an isogeny factor of $V_{\frac{d}{v}}(B/\mathbb{F}_{q^v})$ and~$\phi(\frac{d}{v}) \geq s$.
				\item \label{prop:wr:decomposition:intermediate_field:v:subfield:e}
				For every finite field~$\mathbb{F}_{q^v} \subseteq \mathbb{F}_{q^e}$, the abelian variety~$A/\mathbb{F}_{q^v}$ is an isogeny factor of $W_{e,v} \left(V_{\frac{d}{e}}(B/\mathbb{F}_{q^e})\right)$ and~$s \leq \frac{e}{v}\phi(\frac{d}{e})$.
			\end{enumerate}
		\end{proposition}
		\begin{proof}
			
			\ref{prop:wr:decomposition:intermediate_field:Wde}. The decomposition of~$W_{\mathbb{F}_{q^d}/\mathbb{F}_{q^e}}(B)$ into a product of primitive groups~$V_r(B)$ can be obtained using results from~\S\ref{sec:char_poly_from_prim_groups}.
			
			\ref{prop:wr:decomposition:intermediate_field:Ae}. This follows Proposition~\ref{prop:av:factor:of:Ve}.
			
			\ref{prop:wr:decomposition:intermediate_field:frobA}. From previous statement we have that~$\chi_{A,q^e}(T)$ divides~$\chi_{V_{\frac{d}{e}}(B),q^e}(T)$. From Proposition~\ref{prop:av_is_factor_of_W} it follows that~$\chi_{A,q}(T) \mid \chi_{W_{\mathbb{F}_{q^e}/\mathbb{F}_{q}}(A),q}(T) = \chi_{A,q^e}(T^e)$. Since~$\chi_{A,q^e}(T^e)$ is a divisor of~$\chi_{V_{\frac{d}{e}}(B),q^e}(T^e)$, we obtain that~$\chi_{A,q}(T) \mid \chi_{V_{\frac{d}{e}}(B),q^e}(T^e)$.
			
			\ref{prop:wr:decomposition:intermediate_field:A}. We have~$\chi_{A,q}(T) \mid \chi_{V_{\frac{d}{e}}(B),q^e}(T^e) = \chi_{W_{\mathbb{F}_{q^e}/\mathbb{F}_q}(V_{\frac{d}{e}}(B)),q}(T)$. Therefore, the abelian variety~$A$ is isogenous to an abelian subvariety of~$W_{\mathbb{F}_{q^e}/\mathbb{F}_q}(V_{\frac{d}{e}}(B))$ by Tate's theorem~\cite{Tate1966}.
			
			\ref{prop:wr:decomposition:intermediate_field:Ae:isogV}. Let~$A/\mathbb{F}_{q^e} \sim V_{\frac{d}{e}}(B)$. Then~$g = \dim A = \dim V_{\frac{d}{e}}(B) = \phi(\frac{d}{e}) g / s$ and so~$s = \phi(\frac{d}{e})$.
			
			\ref{prop:wr:decomposition:intermediate_field:Ae:subvarV}. Let~$A/\mathbb{F}_{q^e}$ be isogenous to a proper subvariety of~$V_{\frac{d}{e}}(B)$. Then~$g = \dim A < \dim V_{\frac{d}{e}}(B) = \phi(\frac{d}{e})g/s$ and~$s < \phi(\frac{d}{e})$.
			
			\ref{prop:wr:decomposition:intermediate_field:factor_of_Vte}. This is Proposition~\ref{prop:av:factor:of:Ve} with~$\mathbb{F}_{q^e}$ instead of~$\mathbb{F}_q$.
			
			\ref{prop:wr:decomposition:intermediate_field:v:subfield:e}. We know that $A/\mathbb{F}_{q^v} \mid W_{d,v}(B)$. From properties of the Weil restriction we have $W_{d,v}(B) = W_{e,v}(W_{d,e}(B)) \sim W_{e,v}(W_{d,e}(B)) \sim \prod_{r \mid \frac{d}{e}} W_{e,v}(V_r(B))$. Therefore, $A/\mathbb{F}_{q^v} \mid W_{e,v}(V_r(B))$ for some~$r \mid \frac{d}{e}$. Since $V_r(B) \sim B^{\phi(r)}$ and~$W_{e,v}(V_r(B)) \sim B^{\frac{e}{v} \phi(r)}$ over~$\mathbb{F}_{q^{e r}}$, abelian variety~$A$ is not simple over~$\mathbb{F}_{q^{e r}}$. This is only possible when~$r = \frac{d}{e}$ due to minimality of~$d$. Condition~$s \leq \frac{e}{v}\phi(\frac{d}{e})$ can be obtained from inequality~$\dim A \leq \dim W_{e,v}(V_{\frac{d}{e}}(B))$.
		\end{proof}

		Note that there may be several minimal fields~$\mathbb{F}_{q^e}$. The proposition above allows us to determine~$A$ as an isogeny factor of one of the primitive groups.
		
		\begin{proposition}
			Under the assumptions of Proposition~\ref{prop:wr:decomposition:intermediate_field} let~$4 \nmid e$, and~$A$ be ordinary. Then $W_{\mathbb{F}_{q^d}/\mathbb{F}_q}(B) \sim W_{\mathbb{F}_{q^e}/\mathbb{F}_q}(\widetilde{B}) \times \widetilde{W}$, where~$W_{\mathbb{F}_{q^e}/\mathbb{F}_q}(\widetilde{B})$ is simple and $\widetilde{W}$ is an abelian variety. Moreover, there are only the following three possibilities for~$A$.
			\begin{enumerate}
				\item $A \sim W_{\mathbb{F}_{q^e}/\mathbb{F}_q}(\widetilde{B})$ and $e = s$, $d > s + 1$.
				\item $A$ is a proper factor of $\widetilde{W}$ and $d > s + e + 1$.
				\item $A \sim \widetilde{W}$ and~$d = s + e$.
			\end{enumerate}
		\end{proposition}
		\begin{proof}
			Since~$\widetilde{B}/\mathbb{F}_{q^e}$ is primitive, ordinary, and~$4 \nmid e$, the variety~$W_{\mathbb{F}_{q^e}/\mathbb{F}_q}(\widetilde{B})$ is simple by Proposition~\ref{prop:wr:simplicity}. It is known (see~\cite[\S2.2]{DiemNaumann2003}) that $W_{\mathbb{F}_{q^e}/\mathbb{F}_q}(\widetilde{B})$ is isogenous to a subvariety of~$W_{\mathbb{F}_{q^d}/\mathbb{F}_q}(\widetilde{B}) \sim W_{\mathbb{F}_{q^d}/\mathbb{F}_q}(B)$. Thus, we have $W_{\mathbb{F}_{q^d}/\mathbb{F}_q}(B) \sim W_{\mathbb{F}_{q^e}/\mathbb{F}_q}(\widetilde{B}) \times \widetilde{W}$ for simple~$W_{\mathbb{F}_{q^e}/\mathbb{F}_q}(\widetilde{B})$ and some abelian variety~$\widetilde{W}$. Thus, the abelian variety~$A$ is either isogenous to~$W_{\mathbb{F}_{q^e}/\mathbb{F}_q}(\widetilde{B})$, or to a factor of~$\widetilde{W}$. The first case is possible only when~$\dim W_{\mathbb{F}_{q^e}/\mathbb{F}_q}(\widetilde{B}) = \dim A$. This means that~$e = s$. Moreover, $W_{\mathbb{F}_{q^d}/\mathbb{F}_q}(B) \sim B^d$ over~$\mathbb{F}_{q^d}$, so~$\widetilde{W}$ can not have dimension less than or equal to~$\dim B$. Otherwise, either $B$ is not simple over~$\mathbb{F}_{q^d}$, or it is not primitive over~$\mathbb{F}_{q^e}$. Therefore, $\dim \widetilde{W} = \frac{d g}{s} - \frac{e g}{s} > \dim B = \frac{g}{s}$ and so $d > s + 1$ when~$e = s$. 
			
			Let~$A$ be a proper factor of~$\widetilde{W}$. Then
			\[
			W_{\mathbb{F}_{q^d}/\mathbb{F}_q}(B)
			\sim
			W_{\mathbb{F}_{q^e}/\mathbb{F}_q}(B)
			\times
			A
			\times
			\widetilde{W}',
			\]
			where~$\widetilde{W}'$ is an abelian variety with~$\dim \widetilde{W}' > \dim B$. Therefore, it holds that~$\frac{d g}{s} - \frac{e g}{s} - g > \frac{g}{s}$ and so~$d > s + e + 1$.
			
			Let~$A \sim \widetilde{W}$. Then $g = \dim A = \dim \widetilde{W} = \frac{d g}{s} - \frac{e g}{s}$ and so $d = s + e$.
		\end{proof}
		The proposition above allows us to find the characteristic polynomial~$\chi_{A,q}(T)$ either as~$\chi_{A,q}(T) = \chi_{W_{\mathbb{F}_{q^e}/\mathbb{F}_q}(\widetilde{B}),q}(T) = \chi_{\widetilde{B},q^e}(T^e)$, or as a divisor of
		\[\chi_{\widetilde{W},q}(T)
		=
		\chi_{W_{\mathbb{F}_{q^d}/\mathbb{F}_q}(B),q}(T)
		/
		\chi_{W_{\mathbb{F}_{q^e}/\mathbb{F}_q}(\widetilde{B}),q}(T)
		=
		\chi_{B,q^d}(T)
		/
		\chi_{\widetilde{B},q^e}(T^e)
		.
		\]
		Here possible variants for the characteristic polynomial of~$\chi_{\widetilde{B},q^e}(T)$ can be determined from Proposition~\ref{prop:wr:decomposition:intermediate_field}. Since~$\widetilde{B}$ is a factor of~$W_{\mathbb{F}_{q^d}/\mathbb{F}_{q^e}}(B)$, the abelian variety~$\widetilde{B}$ is isogenous to a subvariety of~$V_{\widetilde{e}}(B)$.

		\subsubsection{The Weil restriction components enumeration}
		In this section, we describe an algorithm that enumerates all possible components of the Weil restriction~$W_{d,1}(B)$ that can contain (as an isogeny factor) a geometrically split abelian variety~$A$ of dimension~$g$. We present it as Algorithm~\ref{alg:wr:components}. The structure of the Weil restriction depends on the minimal field of definition~$\mathbb{F}_{q^e}$ of the isogeny class of~$B$. The algorithm enumerates all possible fields of definition of~$B$ and checks that for every proper subfield~$\mathbb{F}_{q^v}$ of~$\mathbb{F}_{q^d}$ there is a component of the Weil restriction~$W_{d,v}(B)$ of suitable dimension (it must have the dimension greater or equal to the dimension of~$A$). More precisely, if $A$ is a geometrically split abelian variety such that $A \sim B^s$ over~$\mathbb{F}_{q^d}$, $A$ is simple over all proper subfields of~$\mathbb{F}_{q^d}$, and~$B$ has minimal (up to isogeny) field of definition~$\mathbb{F}_{q^e} \subseteq \mathbb{F}_{q^d}$ then $A$ must satisfy the following conditions:
		\begin{enumerate}
			\item $A/\mathbb{F}_{q^v}$ is an isogeny factor~$W_{d,v}(B)$ for all~$\mathbb{F}_{q^v} \subseteq \mathbb{F}_{q^d}$. The algorithm checks that decomposition of~$W_{d,v}(B)$ contains an isogeny factor of dimension~$\geq \dim A$.
			\item The structure of~$W_{d,v}(B)$ is determined by the field of definition of~$B$. This structure is described by Proposition~\ref{prop:av:factor:of:Ve}, Corollary~\ref{cor:wr:primitive:4nd:case}, and Proposition~\ref{prop:wr:primitive:4dcase}.
		\end{enumerate}
		The algorithm checks these conditions and returns a list of possible variants for~$A$ over the field~$\mathbb{F}_{q}$. 
		
		\begin{algorithm}[H]
			\caption{(WResComponents) The Weil restriction components.}\label{alg:wr:components}
			\KwData{A triple~$(g,d,s)$ such that~$A \sim B^s$ over~$\mathbb{F}_{q^d}$, where $A$ is an abelian variety simple over all proper subfields of~$\mathbb{F}_{q^d}$.}
			\KwResult{
				A list of components of the Weil restriction~$W_{d,1}(B)$ that can contain~$A/\mathbb{F}_{q}$ as an isogeny factor.
			}
			Construct a list~$L$ of all binary vectors such that for~$\ell \in L$ and~$v \mid d$:
			if~$\ell_v = 1$ then $ \ell_u = 1$ for all $u \mid d, v \mid u$; if $\ell_v = 0$ then $\forall u \mid v$ it holds that~$\ell_u = 0$\;
			$S = []$\;
			\ForAll{$\ell \in L$}{
				$G = \{\}$\tcp*[f]{key/value dictionary}\;
				\ForAll{$v \mid d$ and $v \neq d$}{
					\uIf{$\ell_v = 1$ and $\phi(\frac{d}{v}) \geq s$}{
						$G[v] = V_{\frac{d}{v}}(B/\mathbb{F}_{q^v})$\;
					}
					\ElseIf{$\ell_v = 0$}{
						find the minimal $u$ such that~$\ell_u = 1$ and~$v \mid u \mid d$\;
							\If{$u = d$ and ($4 \nmid \frac{d}{v}$ and $u = v s$ or $4 \mid \frac{d}{v}$ and ($u \geq  2 v s$ or $u=v s$))\label{alg:wr:components:cond:B:undef}
							}{
								$G[v] = W_{u,v}(B/\mathbb{F}_{q^u})$\;
							}
							\If{$u \neq d$}{
								$G[v] = W_{u,v}(V_{\frac{d}{u}}(B/\mathbb{F}_{q^u}))$\;
							}
						}
						\If{$G[v] = \emptyset$
							\label{alg:wr:components:cond:ell:exclusion}
						}{
							$G = \{\}$; break\;
						}
					}
					\If{$G \neq \{\}$}{
						$S = \operatorname{concat}(S, [G])$\;
					}
				}
				\Return{$\{ G[1] \mid G \in S \}$}
			\end{algorithm}
			
			\begin{proposition}
				\label{prop:wcomponents:correctness}
				Let~$A$ be an abelian variety of dimension~$g$ over a finite field~$\mathbb{F}_{q}$, let~$A \sim B^s$ over~$\mathbb{F}_{q^d}$, and let $A$ be simple over all proper subfields of~$\mathbb{F}_{q^d}$. Given the integers $g \geq 2$, $d \geq 2$, and $s \geq 2$, Algorithm~\ref{alg:wr:components} returns a set of all possible components~$V$ of $W_{d,1}(B)$ such that the following conditions hold.
				\begin{enumerate}
					\item $V = W_{u,1}(V_{\frac{d}{u}}(B/\mathbb{F}_{q^u}))$, where~$\mathbb{F}_{q^u}$ is the minimal field of definition of the isogeny class of~$B$.
					\item $\dim V \geq \dim A$.
					\item For each choice of $V$ it also holds that for every proper subfield $\mathbb{F}_{q^v}$ of $\mathbb{F}_{q^d}$ the Weil restriction~$W_{d,v}(B)$ has a component $W$ such that $\dim W \geq \dim A$.
				\end{enumerate}
			\end{proposition}
			\begin{proof}
				The list $L$ in the algorithm represents all possibilities for compatible fields of definitions of~$B$. For~$\ell \in L$ we put~$\ell_v = 1$ if the isogeny class of~$B$ is defined over~$\mathbb{F}_{q^v}$ and~$\ell_v = 0$ otherwise. If~$B$ is defined over~$\mathbb{F}_{q^v}$ then it is defined over all extensions of~$\mathbb{F}_{q^v}$, so we put also $\ell_u = 1$ for all~$v \mid u \mid d$ when~$\ell_v = 1$. Similarly, if $B$ is undefined over~$\mathbb{F}_{q^v}$ it is undefined over all subfields of~$\mathbb{F}_{q^v}$. So, we put~$\ell_u = 0$ for all $u \mid v$ when~$\ell_v = 0$. Thus, the element $\ell$ contains a complete information about fields of definition of~$B$ over all proper subfields of~$\mathbb{F}_{q^d}$. When we fix some~$\ell$, the decomposition of the Weil restriction~$W_{d,v}(B)$ is determined uniquely by Proposition~\ref{prop:av:factor:of:Ve}, Corollary~\ref{cor:wr:primitive:4nd:case}, and Proposition~\ref{prop:wr:primitive:4dcase}. These three statements give us the component~$V$ of required form. We note that the form $V = W_{u,1}(V_{\frac{d}{u}}(B/\mathbb{F}_{q^u}))$ includes other two variants in the algorithm, since $V_1(C) = C$ and~$W_{1,1}(C) = C$ for any abelian variety~$C$. In general, over a subfield~$\mathbb{F}_{q^v}$ of $\mathbb{F}_{q^d}$, we have a component $W$ of the form $W  = W_{u,v}(V_{\frac{d}{u}}(B/\mathbb{F}_{q^u}))$. So, in the algorithm we enumerate all elements~$\ell \in L$ and check that for all~$v \mid d$ there is a component~$W$ of~$W_{d,v}(B)$ of dimension~$\geq g = \dim A$. If for some~$v$ there is no such a component then we have a contradiction and we exclude such an~$\ell$. In the end we obtain a list of all possible components~$V$ of~$W_{d,1}(B)$, such that~$A \mid V$ and $A$ can also be a factor of some component of~$W_{d,v}(B)$ for all subfields of~$\mathbb{F}_{q^d}$. The component $W$ has dimension~$\dim W = \dim W_{u,v}(V_{\frac{d}{u}}(B/\mathbb{F}_{q^u})) = \frac{u}{v} \cdot \dim V_{\frac{d}{u}}(B/\mathbb{F}_{q^u})) = \frac{u}{v} \cdot \phi(\frac{d}{u}) \cdot \dim B = \frac{u}{v} \cdot \phi(\frac{d}{u}) \cdot \frac{g}{s}$. So, the condition $\dim W \geq g$ is equivalent to~$u \phi(\frac{d}{u}) \geq v s$. At step~\ref{alg:wr:components:cond:B:undef} we have $B$ primitive over $\mathbb{F}_{q^d}$, so $W_{d,v}(B)$ is simple when $4 \nmid \frac{d}{v}$ (by Corollary~\ref{cor:wr:primitive:4nd:case}) or $4 \mid \frac{d}{v}$ and~$\pi_B \not\in -4\mathbb{Q}[\pi_B]^4$. In this case~$A/\mathbb{F}_{q^v} \sim W$ and the dimension check condition is~$u = vs$. When $4 \nmid \frac{d}{v}$ and~$\pi_B \in -4\mathbb{Q}[\pi_B]^4$ we have that~$W_{d,v}(B) \sim W_0 \times \widetilde{W}_0$. Therefore, the dimension check condition is~$\dim W_0 \geq g$ or equivalently $u = d \geq 2 v s$.
				
				At step~\ref{alg:wr:components:cond:ell:exclusion} we exclude the element~$\ell$ if for the field~$\mathbb{F}_{q^v}$ there is no suitable component of~$W_{d,v}(B)$ that can contain $A$ as isogeny factor. Thus, the output of the algorithm contains only factors of $W_{d,1}(B)$ for which there is an element~$\ell$ with existing decompositions over all corresponding subfields~$\mathbb{F}_{q^v}$.
			\end{proof}

			\section{General algorithm for enumeration}
			\label{sec:algorithm}
			In this section, we describe an algorithm that enumerates all possible variants for geometrically split ordinary abelian varieties of dimension~$g$. There are two cases. The first case when~$\chi_{A,q}(T) \in \mathbb{Z}[T^d]$ and~$d$ is the minimal splitting degree, is straightforward and is described in the following theorem.
			\begin{theorem}
				\label{th:av:gsplit:case1} Let~$A$ be a simple ordinary geometrically split abelian variety of dimension~$g$ over a finite field~$\mathbb{F}_{q}$. Assume that~$A \sim B^s$ over~$\mathbb{F}_{q^d}$ for a simple abelian variety~$B$ primitive over~$\mathbb{F}_{q^e} \subseteq  \mathbb{F}_{q^d}$ and~$A$ is simple over all proper subfields of~$\mathbb{F}_{q^d}$. If~$\chi_{A,q}(T) \in \mathbb{Z}[T^d]$ then
				\begin{enumerate}
					\item $d$ is a prime divisor of~$g$;
					\item $e = d$ and~$s = d$;
					\item $A \sim W_{d,1}(B)$;
					\item $\chi_{A,q}(T) = \chi_{B,q^d}(T^d)$.
				\end{enumerate} 
			\end{theorem}
			\begin{proof}
				This follows from Proposition~\ref{prop:splitting_degrees:minimal} except for the fact that~$e=d$. Since $d$ is prime, the only proper subfield of~$\mathbb{F}_{q^d}$ is~$\mathbb{F}_{q}$. If~$e = 1$ then $A$ must be an isogeny factor of~$W_{d,1}(B) \sim~V_{1}(B) \times V_d(B)$. We have~$\dim V_d(B) = \phi(d) \cdot \frac{g}{s} = (d-1) \frac{g}{d} = g - g/d < g$. Therefore, $W_{d,1}(B)$ does not contain a factor of dimension~$g$ when~$B$ is defined over~$\mathbb{F}_q$. Thus, the case~$e = 1$ is impossible and the only remaining case is~$e = d$.
			\end{proof}

			For the second case, we present Algorithm~\ref{alg:geom:split:enum}. This algorithm consists of two main parts. The first is the enumeration of possible endomorphism rings of~$B$ by inspecting the intersection of fields~$\mathbb{Q}(\zeta_d)$ and~$\mathbb{Q}(\pi_B)$. If we have available endomorphism rings for~$B$ then we enumerate and return the Weil restriction components using Algorithm~\ref{alg:wr:components}.
			
			\begin{algorithm}[H]
				\caption{List of simple geometrically split ordinary abelian varieties when~$\mathbb{Q}(\pi) = \mathbb{Q}(\pi^d, \zeta_d)$.}\label{alg:geom:split:enum}
				\KwData{A dimension $g \geq 2$.%
				}
				\KwResult{List of all possible tuples~$(V, d, s, e, \chi_{V,q})$ such that for every ordinary simple geometrically split abelian variety~$A$ there is an element of this list satisfying the following properties.
					\begin{itemize}
						\item $A \sim B^s$ over~$\mathbb{F}_{q^d}$, isogeny class of~$B$ is primitive over~$\mathbb{F}_{q^e}$, and~$A$ is simple over all proper subfields of $\mathbb{F}_{q^d}$
						\item $V$ is a component of $W_{d,1}(B)$ s.t. $A \mid V$;
						\item $\chi_{V,{q}}(T) \in \mathbb{Z}[T, q, b_{1,e}, \ldots b_{\frac{g}{s},e}]$ is the characteristic polynomial of the Frobenius endomorphism on $V/\mathbb{F}_{q^e}$ where $b_{1,e}, \ldots b_{\frac{g}{s},e}$ are coefficients of~$\chi_{B,q}(T)$.
					\end{itemize}
				}
				$L = []$\;
				\For{$(d,s)$ s.~t. $d > 2$, $s \leq \phi(d) \leq 2 g$,  $\phi(d) \mid 2 g$, and~$s \mid \phi(d)$}{
					\If{$\phi(d) > s$ and $\operatorname{EnumEnd(g,s,d)} = \emptyset$}{
						choose next $(d,s)$\;
					}
					$S = \operatorname{WResComponents}(g, d, s)$\;
					\For{$V \in S$}{
						$e = \{\text{field of definition of}~B\}$\;
						$L = \operatorname{concat}(L, [(V, d, s, e, \chi_{V,q})])$\;
					}
				}
				\Return{L}
			\end{algorithm}
			
			\begin{theorem}
				Algorithm~\ref{alg:geom:split:enum}, given the dimension~$g$, returns a list of all possible tuples~$(V, d, s, e, \chi_{V,q})$ such that
				\begin{enumerate}
					\item
					\label{th:alg:general:d_n_s}
					$e \geq 1$, $d>2$, $2 \leq s \leq \phi(d) \leq 2 g$,  $\phi(d) \mid 2 g$, and~$s \mid \phi(d)$.
					\item $V = W_{e,1}(V_{\frac{d}{e}}(B))$ for some simple abelian variety~$B$ of dimension~$\frac{g}{s}$ from primitive~over~$\mathbb{F}_{q^e}$ isogeny class.
					\item\label{th:alg:general:wres:cond1}
					$V$ is a component of $W_{d,1}(B)$ such that~$\dim V \geq g$.
					\item\label{th:alg:general:wres:cond2}
					For each subfield~$\mathbb{F}_{q^v}$ of $\mathbb{F}_{q^d}$ there is a component~$W$ of~$W_{d,v}(B)$ such that~$\dim W \geq g$.
					\item\label{th:alg:general:end:cond1}
					If $\phi(d) = 2 g$ then there is a CM-subfield of $\mathbb{Q}(\zeta_d)$ isomorphic of degree~$\frac{2 g}{s}$.
					\item\label{th:alg:general:end:cond2}
					If $\phi(d) > s$, then there is a subfield~$M$ of~$\mathbb{Q}(\zeta_d)$ such that~$\deg M = \phi(d) / s$. Moreover, if $M$ is real then~$2 \deg M \mid \frac{2 g}{s}$. 
					\item $\chi_{V,q}$ is the characteristic polynomial of~$V$ in terms of coefficients of~$\chi_{B,q}$.
				\end{enumerate}
				
				Moreover, for every ordinary simple geometrically split abelian variety~$A$ such that~$\mathbb{Q}(\pi_A) = \mathbb{Q}(\pi_A^{d_0}, \zeta_{d_0})$ there is a tuple from the list with the following properties.
				\begin{itemize}
					\item $d = d_0$, $A \sim B^s$ over~$\mathbb{F}_{q^d}$.
					\item $A$ is simple over all proper subfields of~$\mathbb{F}_{q^d}$.
					\item $A$ is an isogeny factor of~$V$.
					
				\end{itemize}
			\end{theorem}
			\begin{proof}
				All restrictions on the elements of the list are necessary conditions for ordinary simple geometrically abelian varieties. Lets consider them in more detail. 
				
				The conditions on~$(d,s)$ in the for-loop are taken from Proposition~\ref{prop:spliting_degrees}. Thus, condition~\ref{th:alg:general:d_n_s} of the theorem is the necessary condition for a minimal splitting degree~$d$ and exponent~$s$ of ordinary simple geometrically split abelian variety~$A$ for the case when~$\mathbb{Q}(\pi_A) = \mathbb{Q}(\pi_A^d, \zeta_d)$.
				
				The degree of intersection of $\mathbb{Q}(\pi^d) = \mathbb{Q}(\pi_B)$ and~$\mathbb{Q}(\zeta_d)$ is equal to~$\frac{\phi(d)}{s}$. Thus, the intersection is nontrivial when $\phi(d) > s$. In this case we enumerate compatible endomorphism subrings of~$\mathbb{Q}(\pi_B)$ using Algorithm~\ref{alg:end:enum}. By Proposition~\ref{prop:enum:end:props} after this check we have a triple $(g,s,d)$ satisfying conditions~\ref{th:alg:general:end:cond1} and~\ref{th:alg:general:end:cond2} of our theorem.
				
				If there are suitable endomorphism subrings, then we enumerate components of the Weil restriction~$W_{d,1}(B)$ using Algorithm~\ref{alg:wr:components}. By Proposition~\ref{prop:wcomponents:correctness} it returns a list of all components satisfying conditions~\ref{th:alg:general:wres:cond1} and~\ref{th:alg:general:wres:cond2}. Each such component is of the form~$W_{e,1}(V_{\frac{d}{e}}(B))$, so we can compute $\chi_{V,q}$ by computing $\chi_{V_{\frac{d}{e}}(B), q^e}(T^{\frac{d}{e}})$ using Proposition~\ref{prop:Ve:charpoly}
			\end{proof}

			It is important to note that we provide only necessary conditions and do not claim existence of every abelian variety in the list. 
			
			\subsection{Implementation}
			We implemented\footnote{Source code can be found at~\url{https://github.com/novoselov-sa/avGsplit}} Algorithm~\ref{alg:geom:split:enum} in the SageMath 10.6 computer algebra system. For fast computation of the inverse $\phi$ function we used the invphi package~\cite{Alekseyev2016} for Pari/GP. The enumeration results of geometrically split abelian varieties for dimensions $6-10$ are presented in~\ref{sec:dim6-10}. We consider the cases of dimensions~$3$ and~$4$, as well as prime dimensions, in more detail in the following sections. Current implementation allows us to perform enumeration for dimension up to~$100$ ($73$ cases of geometrically split abelian varieties) and higher. 
			
			\section{Dimension $3$}
			\label{sec:dim3}
			Over the base field an abelian threefold~$A$ can be either simple, or it splits in two ways: $A \sim E \times B$ or~$A \sim E_1 \times E_2 \times E_3$, where~$E, E_1, E_2, E_3$ are elliptic curves and~$B$ is a simple abelian surface. Explicit examples of such a split abelian variety can be found in~\cite{NovBol2022}. In our work we study the case of ordinary geometrically split abelian threefolds and do not consider absolutely simple abelian varieties. The characteristic polynomials of simple ordinary geometrically split abelian varieties of dimension~$3$ are described in the following theorem.
			\begin{theorem}
				Let~$A$ be a simple ordinary geometrically split abelian variety of dimension~$3$ over a finite field~$\mathbb{F}_q$. Assume that $A$ is simple over all proper subfields of~$\mathbb{F}_{q^d}$. Then~$A \sim E^3$ over~$\mathbb{F}_{q^d}$ and the possible values of~$d$ are~$3$ and~$7$. Moreover, let~$t_k$ be such that~$\chi_{E,q^k}(T) = T^2 - t_k T + q^k$ and~$t = t_1$. Then the following holds.
				\begin{enumerate}
					\item If $d=3$ then~$\chi_{A,q}(T) = T^6 - t_3 T^3 + q^3$.
					\item
					\label{th:gsplit:dim3:d7}
					If $d=7$ then~$E$ is an elliptic curve with~$j$-invariant~$-3375$ and~$\chi_{A,q}(T)$ is one of two irreducible factors of the polynomial:
					\begin{equation*}
						\begin{split}
							f_7 = T^{12} &+ v_1 T^{11} + v_2 T^{10} + v_3 T^9 + v_4 T^8 + v_5 T^7+ v_6 T^6+v_5 q T^5 +\\&+ v_4 q^2 T^4 + v_3 q^3 T^3 + v_2 q^4 T^2+q^5 v_1 T+q^6,
						\end{split}
					\end{equation*}
					where~$v_1 = t$, $v_2 = (t^2-q)$, $v_3 = t^3-2 q t$, $v_4 = t^4-3 q t^2+q^2$, $v_5 = t^5-4 q t^3+3 q^2 t$, $v_6 = t^6-5 q t^4+6 q^2 t^2-q^3$.
				\end{enumerate}
			\end{theorem}
			\textit{Proof}.
			By Proposition~\ref{prop:av:isotypicity}, the abelian variety~$A$ can only be decomposed over~$\mathbb{F}_{q^d}$ as~$A \sim E^3$. 
			Consider now the possible degrees~$d$. Let~$\pi$ be the Frobenius endomorphism on~$A$ over~$\mathbb{F}_q$. Then~$\pi^d$ is the Frobenius endomorphism on~$A$ over~$\mathbb{F}_{q^d}$. By Proposition~\ref{prop:spliting_degrees} there are two cases:
			\begin{enumerate}
				\item $\chi_{A,{q}} \in \mathbb{Z}[\pi]$ and~$d = 3$.
				\item $\mathbb{Q}(\pi) = \mathbb{Q}(\pi^d, \zeta_d)$ and~$d = 7, 9, 14, 18$.
			\end{enumerate}
			We consider the two cases separately.
			
			\textit{Case 1}. For~$d = 3$ the characteristic polynomial of~$\pi$ has the form~$\chi_{A,{q}}(T) = T^6 + a_3 T^3 + q^3$. Then it can be shown from known recurrence formulae~\cite[\S17.1.2]{CohFre+2005} for coefficients of~$\chi_{A,q}(T)$ that~$a_3 = -t_3$, where~$t_3$ is a trace of Frobenius for elliptic curve~$E$ over~$\mathbb{F}_{q^3}$.
			
			\textit{Case 2}. Let~$\mathbb{Q}(\pi) = \mathbb{Q}(\pi^d, \zeta_d)$. For~$d = 7, 9, 14, 18$ we have $\phi(d) = 6$ and so $\mathbb{Q}(\pi) = \mathbb{Q}(\zeta_d)$. For the cases~$d = 9, 14$, and~$18$ it can be shown by applying Algorithm~\ref{alg:wr:components} that over a proper subfield of~$\mathbb{F}_{q^d}$ there is no Weil component that have~$A$ as an isogeny factor (by checking dimensions). Since~$A/\mathbb{F}_{q^d} \sim E^3$, the field~$\mathbb{Q}(\pi^d)$ is a quadratic subfield of~$\mathbb{Q}(\pi) = \mathbb{Q}(\zeta_d)$ isomorphic to endomorphism ring of the elliptic curve~$E$. Cyclotomic field $Q(\zeta_7)$ %
			contains only one quadratic subfield -- $Q(\sqrt{-7})$. This corresponds to elliptic curve with $j$-invariant~$-3375$ (this is the only factor of Hilbert modular polynomial for discriminant~$d=-7$).%

			Let~$d = 7$ then~$W \sim V_1(E) \oplus V_7(E)$ and~$A$ is a factor of~$V_7(E)$. The elliptic curve~$E$ is ordinary and its endomorphism ring is a subfield of~$\mathbb{Q}(\zeta_7)$. Then~$V_7(E)$ can be decomposed into two simple non-isogenous abelian varieties of dimension~$3$ by Proposition~3.10 from~\cite{FreemanSatoh2011}. Therefore, our abelian variety~$A$ must be one of them. Thus, the characteristic polynomial~$\chi_{V_7(E),q}(T)$ factors into two distinct irreducible polynomials and one of these polynomials is equal to~$\chi_{A,q}(T)$. By applying the Proposition~\ref{prop:Ve:charpoly}, we obtain that~$\chi_{V_7(E),q}(T) = f_7$. This proves~\ref{th:gsplit:dim3:d7}) of the theorem.

			\qed
			
			The theorem makes it easy to determine the characteristic polynomials of the Frobenius endomorphism when~$d = 7$. We have explicit equations for the curve~$E$ and any twist of~$E$ is suitable for the determining of~$\chi_{A,q}(T)$. An elliptic curve with~$j$-invariant~$-3375$ has equation $y^2=x^3-\frac{125 x}{28}+\frac{125}{28}$ for $p \neq 2,7$. %
			Note that, for~$q=p$ the curve~$E$ with $j$-invariant~$-3375$ can be ordinary only for~$p \equiv 0,1,2,4 \pmod 7$, otherwise it is supersingular (see~\cite[\S2.2]{ACLL+2023}). %
			An example of an abelian variety for~$d=7$ can be found in~\cite[\href{https://www.lmfdb.org/Variety/Abelian/Fq/3/23/al_dp_asj}{3.23.al\_dp\_asj}]{LMFDB}.

			\section{Dimension $4$}
			\label{sec:dim4}
			In dimension~$4$, there are~$9$ types of geometrically split ordinary abelian varieties. They are described in the following theorem.
			
			\begin{theorem}
				\label{th:g4:frobs}
				Let~$A$ be a simple ordinary geometrically split abelian variety of dimension~$4$ defined over a finite field~$\mathbb{F}_q$. Let~$A$ be isogenous to~$B^s$ over~$\mathbb{F}_{q^d}$ and let $A$ be simple over all proper subfields of~$\mathbb{F}_{q^d}$. Assume that~$\mathbb{F}_{q^e} \subseteq \mathbb{F}_{q^d}$ is the minimal extension of~$\mathbb{F}_q$, where the isogeny class of~$B$ is defined and simple over~$\mathbb{F}_{q^d}$. Then all possible~$A$ are given in the Table~\ref{table:av_gsplit_dim4}. The characteristic polynomials of the Frobenius endomorphism for the cases in the table are the following.
				\begin{enumerate}
					\item $d = 2$, $\chi_{A,q}(T) = \chi_{B,q^2}(T^2)$.
					\item $d = 3$, $\chi_{A,q}(T) = \chi_{V_3(B),q}(T)$.
					\item $d = 4$, $\chi_{A,q}(T)$ is a divisor of $\chi_{B,q^4}(T^4)$, $\chi_{A,{q^2}}(T) = \chi_{B,q^4}(T^2)$.
					\item $d = 5$, $\dim B = 1$, $\chi_{A,q}(T) = \chi_{V_5(B),q}(T)$.
					\item $d = 5$, $\dim B = 2$, $\chi_{A,q}(T)$ is a divisor of  $\chi_{V_5(B),q}(T)$. Endomorphism algebra~$\operatorname{End}^0 B$ contains a subfield isomorphic to~$\mathbb{Q}(\sqrt{5})$.
					\item $d = 6$, $\chi_{A,q}(T)$ is a divisor of $\chi_{V_3(B/\mathbb{F}_{q^2}),q^2}(T^2)$, $\chi_{A,q^2}(T) = \chi_{V_3(B/\mathbb{F}_{q^2})}(T)$, $\chi_{A,q^3}(T) = \chi_{B,q^6}(T^2)$. The isogeny class of~$B$ is undefined over~$\mathbb{F}_{q^3}$.
					\item $d = 8$, $\chi_{A,q}(T)$ is a divisor of $\chi_{B,q^8}(T^8)$, $\chi_{A,q^4}(T) = \chi_{B,q^8}(T^2)$.
					Endomorphism algebra~$\operatorname{End}^0 B$ contains a subfield isomorphic to~$\mathbb{Q}(\sqrt{8})$, $\mathbb{Q}(\sqrt{-4})$, or~$\mathbb{Q}(\sqrt{-8})$.
					\item $d = 10$, $\chi_{A,q}(T)$ is a divisor of  $\chi_{V_5(B/\mathbb{F}_{q^2}),q^2}(T^2)$, $\chi_{A,{q^5}}(T) = \chi_{B,q^{10}}(T^2)$. Endomorphism algebra~$\operatorname{End}^0 B$ contains a subfield isomorphic to~$\mathbb{Q}(\sqrt{5})$.
					\item $d = 12$, $\chi_{A,q}(T)$ is a divisor of  $\chi_{V_3(B/\mathbb{F}_{q^4}),q^4}(T^4)$, $\chi_{A,q^4}(T) = \chi_{V_3(B/\mathbb{F}_{q^4}), q^4}(T)$, $\chi_{A,q^6}(T) = \chi_{B,q^{12}}(T^2)$. Endomorphism algebra~$\operatorname{End}^0 B$ contains a subfield isomorphic to~$\mathbb{Q}(\sqrt{-3})$, $\mathbb{Q}(\sqrt{-4})$, or~$\mathbb{Q}(\sqrt{12})$.
				\end{enumerate}
			\end{theorem}
			\begin{proof}
				This result can be obtained by running Algorithm~\ref{alg:geom:split:enum}. For completeness, however, we provide a full proof here. First, let us determine the list of possible degrees~$d$. For geometrically split abelian varieties there are two cases~\cite{HoweZhu2002}:
				\begin{enumerate}
					\item $\chi_{A,q}(T) \in \mathbb{Z}[T^d]$. By Theorem~\ref{th:av:gsplit:case1}, we have~$A \sim W_{2,1}(B)$, $d=2$, $e=2$, and~$\dim B = 2$ in this case.
					\item $\mathbb{Q}(\pi_A) = \mathbb{Q}(\pi_A^d, \zeta_d)$. By Proposition~\ref{prop:spliting_degrees}, we have~$d > 2$ and~$\phi(d) \mid 2 g = 8$. All such degrees~$d$ belong to the set~$\{3, 4, 5, 6, 8, 10, 12\}$. Moreover, by Proposition~\ref{prop:spliting_degrees}, we have~$s \mid \phi(d)$. Therefore, we obtain the following list of possible pairs $(d, s)$: $(3, 2), (4, 2), (5, 2), (5, 4), (6, 2), (8, 2), (8, 4)$, $(10, 2), (10, 4), (12, 2),(12, 4)$.
				\end{enumerate}
				In all cases,~$A$ is an isogeny factor of~$W_{d,1}(B)$. The structure of~$W_{d,1}(B)$ depends on the minimal field of definition of the isogeny class of~$B$. Denote it by~$\mathbb{F}_{q^e} \subset \mathbb{F}_{q^d}$. Then $A/\mathbb{F}_{q^e} \mid V_{\frac{d}{e}}(B/\mathbb{F}_{q^e})$ by Proposition~\ref{prop:av:factor:of:Ve} and thus~$A \mid W_{e,1}(V_{\frac{d}{e}}(B/\mathbb{F}_{q^e}))$. The result now follows by enumerating all possible degrees~$e$.
				
				For~$d = 3$ we have $s=2$ and two cases for the minimal field of definition~$\mathbb{F}_{q^e}$ of the isogeny class of~$B$ -- $\mathbb{F}_{q}$ and~$\mathbb{F}_{q^3}$. For~$e=1$, we have~$A \mid W_{3,1} \sim V_{1}(B) \times V_{3}(B)$ and so~$A \sim V_{3}(B)$. The case~$e=3$ is impossible, since in this case~$W_{3,1}(B)$ is simple by Proposition~\ref{prop:wr:simplicity} and~$\dim W_{3,1}(B) = 3 \cdot \dim B = 6$.
				
				For~$d=4$, we have~$s = 2$ and~$e \in \{1, 2, 4\}$. For~$e=1,2$ we have~$W_{4,2}(B) \sim V_1(B/\mathbb{F}_{q^2}) \times V_2(B/\mathbb{F}_{q^2})$ and $W_{4,2}(B)$ contains no isogeny factor of dimension~$4$, so we exclude this case. In the remaining $e=4$ case the Weil restriction~$W_{4,1}(B)$ can contain an isogeny factor of dimension~$4$ only when~$\pi_B \in -4\mathbb{Q}[\pi_B]^4$ (by Proposition~\ref{prop:wr:simplicity}).
				
				All other cases for~$d$ can be treated similarly by enumerating all possible cases for~$e$ and performing the Weil restriction decomposition and checking the dimensions. Endomorphism algebra data is obtained by enumeration of the subfields of~$\mathbb{Q}(\zeta_d)$ of degree equal to the degree of~$\mathbb{Q}(\zeta_d) \cap \mathbb{Q}(\pi_B)$.
			\end{proof}
			
			\begin{table}[]
				\caption{Simple ordinary geometrically split abelian varieties of dimension~$4$ over~$\mathbb{F}_q$.}
				\label{table:av_gsplit_dim4}
				\begin{tabular}{cccccc}
					\hline
					$A$ & $d$ & $A/\mathbb{F}_{q^d}$ & $e$ & $\dim B$ & conditions \\ \hline
					$W_{2,1}(B)$ & $2$ & $B^2$ & $2$ & $2$ &  \\
					$V_{3}(B)$ & $3$ & $B^2$ & $1$ & $2$ & \\
					$\mid W_{4,1}(B)$ & $4$ & $B^2$ & $4$ & $2$ &  $\pi_B \in -4\mathbb{Q}[\pi_B]^4$, $A/\mathbb{F}_{q^2} \sim W_{4,2}(B)$ \\
					$V_5(B)$ & $5$ & $B^4$ & $1$ & $1$ & \\
					$\mid V_5(B)$ & $5$ & $B^2$ & $1$ & $2$ & \\
					$\mid W_{2,1}(V_3(B))$ & $6$ & $B^2$ & $2$ & $2$ & $A/\mathbb{F}_{q^2} \sim V_3(B/\mathbb{F}_{q^2}), A/\mathbb{F}_{q^3} \sim W_{6,3}(B)$ \\
					$\mid W_{8,1}(B)$ & $8$ & $B^2$ & $8$ & $2$ & $\pi_B \in -4\mathbb{Q}[\pi_B]^4$, $A/\mathbb{F}_{q^4} \sim W_{8,4}(B)$\\
					$\mid W_{2,1}(V_5(B/\mathbb{F}_{q^2}))$ & $10$ & $B^2$ & $2$ & $2$ &  $A/\mathbb{F}_{q^5} \sim W_{10,5}(B)$\\
					$\mid W_{4,1}(V_3(B/\mathbb{F}_{q^4}))$ & $12$ & $B^2$ & $4$ & $2$ &  $A/\mathbb{F}_{q^4} \sim V_3(B/\mathbb{F}_{q^4}), A/\mathbb{F}_{q^6} \sim W_{12,6}(B)$ \\
					\hline
				\end{tabular}
			\end{table}

			The characteristic polynomials in Theorem~\ref{th:g4:frobs} can be written in an explicit form. Suppose that the characteristic polynomial of~$B/\mathbb{F}_{q^e}$ is as follows:
			\begin{itemize}
				\item $\chi_{B,q^e}(T) = T^2 + b_{1,e} T + q^e$ for $\dim B = 1$,
				\item $\chi_{B,q^e}(T) = T^4 + b_{1,e} T^3 + b_{2,e} T^2 + b_{1,e} q^e T + q^{2 e}$ for $\dim B = 2$.
			\end{itemize}
			Then~$\chi_{A,q}(T)$ can be expressed in terms of the coefficients of~$\chi_{B,q^e}(T)$ using Proposition~\ref{prop:Ve:charpoly} for characteristic polynomials of abelian varieties~$V_u(B/\mathbb{F}_{q^e})$:
			\begin{enumerate}
				\item for~$V_3$ and~$\dim{B} = 2$, we have
				\[\chi_{V_3(B/\mathbb{F}_{q^e}),q^e}(T) =
				T^8 - b_{1,e} T^7 + 
				v_2 T^6 +
				v_3 T^5 + 
				v_4 T^4 + 
				v_3 q^e T^3 + 
				v_2 q^{2 e} T^2
				-b_{1,e} q^{3 e} T
				+
				q^{4 e},\] where $v_2 = b_{1,e}^2-b_{2,e}$, $v_3 = 2 b_{1,e} q^e - b_{1,e} b_{2,e}$ and $v_4 = b_{2,e}^2-b_{1,e}^2 q^{e}-q^{2 e}$;
				\item for~$V_5$ when~$\dim{B} = 1$, we obtain
				\[\chi_{V_5(B/\mathbb{F}_{q^e}),q^e}(T) = 
				T^8
				- b_{1,e} T^7 +
				v_2 T^6 +
				v_3 T^5 + 
				v_4 T^4+
				v_3 q^e T^3+
				v_2 q^{2 e} T^2
				-b_{1,e} q^{3 e} T
				+ q^{4 e},\] where $v_2 = b_{1,e}^2-q^e$, $v_3 = 2 b_{1,e} q^e - b_{1,e}^3$ and $v_4 = b_{1,e}^4 - 3 b_{1,e}^2 q^e + q^{2 e}$;
				\item for~$V_5$ and $\dim{B} = 2$, we have\\
				$\chi_{V_5(B/\mathbb{F}_{q^e}),q^e}(T) = 
				T^{16}
				- b_{1,e} 
				T^{15} +
				v_2
				T^{14} + 
				v_3
				T^{13} +
				v_4
				T^{12} + 
				v_5
				T^{11} + 
				v_6
				T^{10} + 
				v_7
				T^9 + 
				v_8
				T^8 + 
				v_7 q^e
				T^7 + 
				v_6 q^{2 e}
				T^6 +
				v_5 q^{3 e}
				T^5 + 
				v_4 q^{4 e}
				T^4 + 
				v_3 q^{5 e}
				T^3+
				v_2 q^{6 e}
				T^2
				-b_{1,e} q^{7 e}
				T+ 
				q^{8 e},
				$
				where
				\begin{align*}
					v_2 &= b_{1,e}^2-b_{2,e}\\
					v_3 &= -b_{1,e}^3+2 b_{1,e} b_{2,e}-b_{1,e} q^e, \\
					v_4 &= b_{1,e}^4 - 3 b_{1,e}^2 b_{2,e} + b_{2,e}^2+2 b_{1,e}^2 q^e-q^{2 e}, \\
					v_5 &= -b_{1,e}^3 b_{2,e} + 2 b_{1,e} b_{2,e}^2 + 2 b_{1,e}^3 q^e-3 b_{1,e} b_{2,e} q^e-3 b_{1,e} q^{2 e},\\
					v_6 &= b_{1,e}^2 b_{2,e}^2 - b_{2,e}^3 - b_{1,e}^4 q^e-b_{1,e}^2 b_{2,e} q^e+3 b_{1,e}^2 q^{2 e}+2 b_{2,e} q^{2 e},\\
					v_7 &= -b_{1,e} b_{2,e}^3+2 b_{1,e}^3 b_{2,e} q^e+2 b_{1,e} b_{2,e}^2 q^e-4 b_{1,e}^3 q^{2 e}-b_{1,e} b_{2,e} q^{2 e}+2 b_{1,e} q^{3 e},\\
					v_8 &= b_{2,e}^4-3 b_{1,e}^2 b_{2,e}^2 q^e+b_{1,e}^4 q^{2 e}+4 b_{1,e}^2 b_{2,e} q^{2 e}-3 b_{2,e}^2 q^{2 e}-b_{1,e}^2 q^{3 e}+q^{4 e}.
				\end{align*}
			\end{enumerate}

			For every case in the theorem, we give an explicit example of an isogeny class of ordinary simple geometrically split abelian varieties (see \ref{sec:examples}).
			
			\section{Prime dimension $g \geq 5$}
			\label{sec:prime:dim}
			When the dimension~$g$ of~$A$ is prime, the classification of geometrically split abelian varieties is simple, because the list of possible splitting degrees~$d$ is short.
			
			\begin{theorem}
				\label{th:av:gsplit:prime:dim}
				Let~$A$ be an ordinary simple geometrically split abelian variety of prime dimension~$g \geq 5$ defined over a finite field~$\mathbb{F}_{q}$. Assume that~$A \sim B^s$ over~$\mathbb{F}_{q^d}$ and that~$A$ is simple over all proper subfields of~$\mathbb{F}_{q^d}$. Then the abelian variety~$A$ satisfies one of the following conditions.
				\begin{enumerate}
					\item $A \sim W_{g,1}(B)$, $A \sim B^g$ over~$\mathbb{F}_{q^g}$, $B$ is primitive over~$\mathbb{F}_{q^g}$, and~$\chi_{A,q}(T) = \chi_{B,q^g}(T^g) = T^{2g} + b_{g,g} T^{g} + q^g$;
					\item $A$ is an isogeny factor of~$V_{2g+1}(B)$ for prime~$2g+1$, $A \sim B^{g}$ over~$\mathbb{F}_{q^{2g + 1}}$, $B$ is defined over~$\mathbb{F}_{q}$, and~$\chi_{A,q}(T)$ is one of two irreducible factors of~$\chi_{V_{2g+1}(B),q}(T)$.
				\end{enumerate}
			\end{theorem}
			\begin{proof}
				Let~$\mathbb{F}_{q^e} \subseteq \mathbb{F}_{q^d}$ be the minimal field of definition of the isogeny class of~$B$. For a geometrically split abelian variety~$A$, we have the following two cases.
				\begin{enumerate}
					\item $\chi_{A,q}(T) \in \mathbb{Z}[T^d]$. In this case~$d = g$, $e=d=g$, $s=d=g$ by Theorem~\ref{th:av:gsplit:case1}. This proves the first part of theorem.
					\item $\mathbb{Q}(\pi_A) = \mathbb{Q}(\pi_A^d, \zeta_d) = \mathbb{Q}(\pi_B, \zeta_d)$. In this case~$\phi(d) \mid 2 g$ and~$s \mid \phi(d)$ by Proposition~\ref{prop:spliting_degrees}. Since~$g$ is prime, we have~$\phi(d) \in \{2, g, 2 g\}$. From~$\phi(d) = 2$ it follows that~$\mathbb{Q}(\pi_A) = \mathbb{Q}(\pi_B, \zeta_d) \simeq \mathbb{Q}(\pi_B)$, so we exclude this case due to the contradiction. The values of Euler function~$\phi$ are even for all~$d \geq 3$, therefore there is no possible~$d$ for~$\phi(d) = g$ when~$g > 2$ is prime. Thus we have~$d$ such that~$\phi(d) = 2 g$. By \cite[Th.~3.1, Prop.~6.1]{Coleman2012} it holds that~$\phi^{-1}(2 g)$ is empty when~$2g+1$ is composite and~$\phi^{-1}(2 g) \in \{2g+1, 4g+2\}$ when~$2g+1$ is prime and~$g \geq 5$. Thus~$d \in \{2g+1, 4g+2\}$. Since~$\phi(d) = 2g$, it holds that~$\mathbb{Q}(\pi_A) = \mathbb{Q}(\zeta_{2g})$ and it follows from Proposition~\ref{prop:splitting_degrees:minimal} that~$s = g$ and so~$\dim B = 1$. For prime~$d = 2 g + 1$, we have~$e \in \{1, 2g+1\}$. When~$e=1$ we have~$A \mid W_{d,1}(B) \sim V_1(B) \times V_d(B)$. Therefore, $A \mid V_{2g+1}(B)$ and~$\dim V_{2g+1}(B) = 2 g$. If $e = 2 g + 1$ then $W_{d,1}(B)$ is simple and can not contain $A$ as an isogeny factor. Consider now $d = 4 g + 2$. We have here~$e \in \{1, 2, 2 g + 1, 4 g + 2\}$. If~$e=1$ then $A/\mathbb{F}_{q^{2g+1}} \mid W_{4 g + 2, 2 g+1}(B) \sim V_{1}(B/\mathbb{F}_{q^{2g+1}}) \times V_{2}(B/\mathbb{F}_{q^{2g+1}})$. Since~$\dim V_2(B) = 1$, the case~$e=1$ is impossible. Let~$e = 2$. If $B$ is not defined over~$\mathbb{F}_{q^{2g+1}}$ then~$W_{4g+2, 2g+1}(B)$ has dimension~$2$ and therefore can not contain~$A$, contradiction. If $B$ is defined over~$\mathbb{F}_{q^{2g+1}}$ then~$W_{4g+2, 2g+1}(B) \sim V_1 \times V_2$ and $A$ is also is not an isogeny factor of~$W_{4g+2, 2g+1}(B)$ over~$\mathbb{F}_{q^{2g+1}}$. The same is true for~$e = 2 g + 1$, so we exclude it. In the remaining~$e = 4 g + 2$ case~$W_{4g+2, 1}(B)$ is simple and so $A \nmid W_{4g+2, 1}(B)$, contradiction. Thus, the only possible case is~$d = 2 g + 1$, $e=1$.
				\end{enumerate}
			\end{proof}
			An explicit example of a geometrically split ordinary abelian variety of dimension~$5$ can be found in~\ref{sec:examples}.	
			
			\section*{Conclusion}
			In this work, we have presented a partial classification of isogeny classes of ordinary simple geometrically split abelian varieties. The following questions remain open for future work.
			\begin{enumerate}
				\item \textit{Decomposing primitive groups $V_v(B)$.} In the case when $\mathbb{Q}(\pi_A) = \mathbb{Q}(\pi^d, \zeta_d)$ and $M = \mathbb{Q}(\pi_B) \cap \mathbb{Q}(\zeta_d) \neq \mathbb{Q}$, we know that the primitive group~$V_v(B)$ can be decomposed in the following way. Let $\Phi_v(T)$ factors as
				\[
				\Phi_v(T) = \prod\limits_{i=1}^{r} \Phi_v^{(i)}(T)
				\]
				over~$M$ for some $r \geq 1$. The field~$M$ can be found explicitly for fixed~$d$ by enumerating subfields (there may be multiple variants for $M$) of $\mathbb{Q}(\zeta_d)$ of degree $[M:\mathbb{Q}] = \frac{\phi(d)}{s}$. Thus, given~$d$ and~$s$, we can factor $\Phi_v(T)$ over a subfield~$M$ of~$\mathbb{Q}(\zeta_d)$. By a result of Diem and Naumann~\cite[\S3.5]{DiemNaumann2003} we know that
				\[
				V_v(B) \sim \prod_{i=1}^{r} V_{v}^{(i)},
				\]
				where~$V_{v}^{(i)}$ is an abelian variety with~$\operatorname{End}^0{V_{v}^{(i)}} \simeq \operatorname{End}^0{B}[T] / (\Phi^{(i)}_v(T))$.
				A geometrically split abelian variety~$A$ such that~$A/\mathbb{F}_{q^d} \sim B^s$ is an isogeny factor of one of these varieties. However, determining the characteristic polynomial of~$V_{v}^{(i)}$ is an open problem.
				\item \textit{Explicit formulae for special families of curves.} There are known families~\cite{Novoselov2020} of hyperelliptic curves with geometrically split Jacobians. The classification from our work can be used to obtain explicit formulae for the Jacobian order. The classification can be combined with considering of twists decompositions as it was done in~\cite{Satoh2009,GuillevicVergnaud2012}.
				\item \textit{Classification of non-ordinary geometrically split abelian varieties.} Theorems and propositions from this work can be generalized to non-ordinary abelian varieties.
			\end{enumerate}
			
			This work is an extended and refined version of preliminary results presented by the author at the SibeCrypt'24 conference~\cite{Novoselov2024} (for dimension~$3$) and the SibeCrypt'26 conference~\cite{Novoselov2026} (for dimension~$4$).

			\section*{Acknowledgment}
			The work is supported by the Ministry of Science and Higher Education of the Russian Federation (Agreement 075-02-2026-1311).
			
			\pagebreak
			\appendix
			
			\section{Geometrically split ordinary abelian varieties of dimensions $6-10$.}
			\label{sec:dim6-10}
			
			The tables in the following sections are outputs of Algorithm~\ref{alg:geom:split:enum}. The case of dimension~$7$ is completely described by Theorem~\ref{th:av:gsplit:prime:dim}, so we did not include it here. The examples were obtained using the methods from~\S\ref{sec:examples}.
			
			\subsection{Dimension $6$}
			\label{sec:dim6}
			
			\begin{table}[H]
				\caption{Simple ordinary geometrically split abelian varieties of dimension 6}
				\label{table:dim6}
				\begin{tabular}{@{}ccccccc@{}}
					\toprule
					$A$                                     & $d$  & $A/\mathbb{F}_{q^d}$ & $e$ & $\dim{B}$ & $\chi_{A,q}(T)$                                & conditions                                                                                                                                             \\ \midrule
					$W_{2,1}(B)$           & $2$  & $B^2$                & 2   & 3         & $\chi_{B,{q^2}}(T^2)$           &                                                                                                                                                        \\
					$W_{3,1}(B)$           & $3$  & $B^3$                & 3   & 2         & $\chi_{B,q^3}(T^3)$           &                                                                                                                                                        \\
					$V_3(B)$                                & $3$  & $B^2$                & 1   & 3         & $\chi_{V_3(B),q}(T)$                           &                                                                                                                                                        \\
					$\mid W_{4,1}(B)$      & $4$  & $B^2$                & 4   & 3         & $\mid \chi_{B,{q^4}}(T^4)$      & $\pi_B \in -4\mathbb{Q}[\pi_B]^4$                                                                                                                      \\
					$\mid W_{2,1}(V_3(B/\mathbb{F}_{q^2}))$ & $6$  & $B^2$                & 2   & 3         & $\mid \chi_{V_3(B/\mathbb{F}_{q^2}),q^2}(T^2)$ & $A/\mathbb{F}_{q^3} \sim W_{6,3}(B)$                                                                                                  \\
					$V_7(B)$                                & $7$  & $B^6$                & 1   & 1         & $\chi_{V_7(B),q}(T)$                           &                                                                                                                                                        \\
					$\mid V_7(B)$                           & $7$  & $B^3$                & 1   & 2         & $\mid \chi_{V_7(B),q}(T)$                      & $B$ is geom. split                                                                                                                                     \\
					$\mid V_7(B)$                           & $7$  & $B^2$                & 1   & 3         & $\mid \chi_{V_7(B),q}(T)$                      &                                                                                                                                                        \\
					$\mid W_{8,1}(B)$      & $8$  & $B^2$                & 8   & 3         & $\mid \chi_{B,q^8}(T^8)$      & \begin{tabular}[c]{@{}c@{}}$\pi_B \in -4\mathbb{Q}[\pi_B]^4$, \\ $A/\mathbb{F}_{q^4} \sim W_{8,4}(B)$\end{tabular}                    \\
					$\mid V_9(B)$                           & $9$  & $B^2$                & 1   & 3         & $\mid \chi_{V_9(B),q}(T)$                      &                                                                                                                                                        \\
					$\mid W_{3,1}(V_3(B/\mathbb{F}_{q^3}))$ & $9$  & $B^2$                & 3   & 3         & $\mid \chi_{V_3(B/\mathbb{F}_{q^3}),q^3}(T^3)$ &                                                                                                                                                        \\
					$\mid W_{4,1}(V_3(B/\mathbb{F}_{q^4}))$ & $12$ & $B^2$                & 4   & 3         & $\mid \chi_{V_3(B/\mathbb{F}_{q^4}),q^4}(T^4)$ & $A/\mathbb{F}_{q^6} \sim W_{12,6}(B)$                                                                                              \\
					$\mid V_{13}(B)$                        & $13$ & $B^3$                & 1   & 2         & $\mid \chi_{V_{13}(B),q}(T)$                     & $B$ is abs. simple                                                                                                                               \\
					$\mid W_{2,1}(V_7(B/\mathbb{F}_{q^2}))$ & $14$ & $B^2$                & 2   & 3         & $\mid \chi_{V_7(B/\mathbb{F}_{q^2}),q^2}(T^2)$ & $A/\mathbb{F}_{q^7} \sim W_{14,7}(B)$                                                                                              \\
					$\mid W_{2,1}(V_9(B/\mathbb{F}_{q^2}))$ & $18$ & $B^2$                & 2   & 3         & $\mid \chi_{V_9(B/\mathbb{F}_{q^2}),q^2}(T^2)$ & \begin{tabular}[c]{@{}c@{}}$A/\mathbb{F}_{q^9} \sim W_{18,9}(B)$,\\ $A/\mathbb{F}_{q^6} \sim V_3(B/\mathbb{F}_{q^6})$\end{tabular} \\
					$\mid W_{6,1}(V_3(B/\mathbb{F}_{q^6}))$ & $18$ & $B^2$                & 6   & 3         & $\mid \chi_{V_3(B/\mathbb{F}_{q^6}),q^6}(T^6)$ & $A/\mathbb{F}_{q^9} \sim W_{18,9}(B)$                                                                                              \\
					$\mid W_{3,1}(V_7(B/\mathbb{F}_{q^3}))$ & $21$ & $B^3$                & 3   & 2         & $\mid \chi_{V_7(B/\mathbb{F}_{q^3}),q^3}(T^3)$ & \begin{tabular}[c]{@{}c@{}}B is geom. split, \\ $A/\mathbb{F}_{q^7} \sim W_{21,7}(B)$
					\end{tabular}                                 \\ \bottomrule
				\end{tabular}
			\end{table}

			\subsection{Dimension $8$}
			\label{sec:dim8}
			\begin{table}[H]
				\caption{Simple ordinary geometrically split abelian varieties of dimension $8$}
				\label{table:dim8}
				\begin{tabular}{@{}lllllll@{}}
					\toprule
					$d$  & $A$                        & $A/\mathbb{F}_{q^d}$ & $e$ & $\dim B$ & $\chi_{A,q}(T)$                     & conditions               \\ \midrule
					$2$  & $W_{2,1}(B)$         & $B^2$       & 2   & 4        & $\chi_{B,{q^2}}(T^2)$         &                          \\
					$3$  & $V_3(B)$                   & $B^2$       & 1   & 4        & $\chi_{V_3(B),q}(T)$              &                          \\
					$4$  & $\mid W_{4,1}(B)$       & $B^2$       & 4   & 4        & $\mid \chi_{B,{q^4}}(T^4)$       & $\pi_B \in -4\mathbb{Q}[\pi_B]^4$ \\
					$5$  & $V_5(B)$                   & $B^4$       & 1   & 2        & $\chi_{V_5(B),q}(T)$              &                          \\
					$5$  & $\mid V_5(B)$                 & $B^2$       & 1   & 4        & $\mid \chi_{V_5(B),q}(T)$            &                          \\
					$6$  & $\mid W_{2,1}(V_3(B/\mathbb{F}_{q^2}))$  & $B^2$       & 2   & 4        & $\mid \chi_{V_3(B/\mathbb{F}_{q^2}),q^2}(T^2)$  &                          \\
					$8$  & $\mid W_{8,1}(B)$       & $B^2$       & 8   & 4        & $\mid \chi_{B,{q^8}}(T^8)$       & $\pi_B \in -4\mathbb{Q}[\pi_B]^4$ \\
					$10$ & $\mid W_{2,1}(V_5(B/\mathbb{F}_{q^2}))$  & $B^2$       & 2   & 4        & $\mid \chi_{V_5(B/\mathbb{F}_{q^2}),q^2}(T^2)$  &                          \\
					$12$ & $\mid W_{4,1}(V_3(B/\mathbb{F}_{q^4}))$  & $B^2$       & 4   & 4        & $\mid \chi_{V_3(B/\mathbb{F}_{q^4}),q^4}(T^4)$  &                          \\
					$15$ & $\mid V_{15}(B)$                & $B^2$       & 1   & 4        & $\mid \chi_{V_{15}(B),q}(T)$           &                          \\
					$15$ & $\mid W_{3,1}(V_5(B/\mathbb{F}_{q^3}))$  & $B^2$       & 3   & 4        & $\mid \chi_{V_5(B/\mathbb{F}_{q^3}),q^3}(T^3)$  &                          \\
					$16$ & $\mid W_{16,1}(B)$     & $B^2$       & 16  & 4        & $\mid \chi_{B,q^{16}}(T^{16})$    & $\pi_B \in -4\mathbb{Q}[\pi_B]^4$ \\
					$20$ & $\mid W_{4,1}(V_5(B/\mathbb{F}_{q^4}))$  & $B^2$       & 4   & 4        & $\mid \chi_{V_5(B/\mathbb{F}_{q^4}),q^4}(T^4)$  &                          \\
					$24$ & $\mid W_{8,1}(V_3(B/\mathbb{F}_{q^8}))$  & $B^2$       & 8   & 4        & $\mid \chi_{V_3(B/\mathbb{F}_{q^8}),q^8}(T^8)$  &                          \\
					$30$ & $\mid W_{2,1}(V_{15}(B/\mathbb{F}_{q^2}))$ & $B^2$       & 2   & 4        & $\mid \chi_{V_{15}(B/\mathbb{F}_{q^2}),q^2}(T^2)$ &                          \\
					$30$ & $\mid W_{6,1}(V_5(B/\mathbb{F}_{q^6}))$  & $B^2$       & 6   & 4        & $\mid \chi_{V_5(B/\mathbb{F}_{q^6}),q^6}(T^6)$  &                          \\ \bottomrule
				\end{tabular}
			\end{table}

			\subsection{Dimension $9$}
			\label{sec:dim9}
			\begin{table}[H]
				\caption{Simple ordinary geometrically split abelian varieties of dimension $9$}
				\label{table:dim9}
				\begin{tabular}{@{}ccccccc@{}}
					\toprule
					$d$  & $A$                & $A/\mathbb{F}_{q^d}$ & $e$ & $\dim B$ & $\chi_{A,q}(T)$           & conditions \\ \midrule
					$3$  & $W_{3,1}(B)$ & $B^3$       & 3   & 3        & $\chi_{B,q^3}(T^3)$ &            \\
					$7$  & $\mid V_7(B)$         & $B^3$       & 1   & 3        & $\mid \chi_{V_7(B),q}(T)$    &            \\
					$19$ & $\mid V_{19}(B)$        & $B^9$       & 1   & 1        & $\mid \chi_{V_{19}(B),q}(T)$   &            \\
					$19$ & $\mid V_{19}(B)$        & $B^3$       & 1   & 3        & $\mid \chi_{V_{19}(B),q}(T)$   &            \\ \bottomrule
				\end{tabular}
			\end{table}

			\subsection{Dimension $10$}
			\label{sec:dim10}
			\begin{table}[H]
				\caption{Simple ordinary geometrically split abelian varieties of dimension $10$}
				\label{table:dim10}
				\begin{tabular}{@{}ccccccc@{}}
					\toprule
					$d$  & $A$                            & $A/\mathbb{F}_{q^d}$ & $e$ & $\dim B$ & $\chi_{A,q}(T)$                       & conditions               \\ \midrule
					$2$  & $W_{2,1}(B)$           & $B^2$       & 2   & 5        & $\chi_{B,q^2}(T^2)$           &                          \\
					$3$  & $V_3(B)$                       & $B^2$       & 1   & 5        & $\chi_{V_3(B),q}(T)$                  &                          \\
					$4$  & $\mid W_{4,1}(B)$         & $B^2$       & 4   & 5        & $\mid \chi_{B,{q^4}}(T^4)$       & $\pi_B \in -4\mathbb{Q}[\pi_B]^4$ \\
					$5$  & $W_{5,1}(B)$             & $B^5$       & 5   & 2        & $\chi_{B,q^5}(T^5)$             &                          \\
					$6$  & $\mid W_{2,1}(V_3(B/\mathbb{F}_{q^2}))$    & $B^2$       & 2   & 5        & $\mid \chi_{V_3(B/\mathbb{F}_{q^2}),q^2}(T^2)$    &                          \\
					$8$  & $\mid W_{8,1}(B)$           & $B^2$       & 8   & 5        & $\mid \chi_{B,q^8}(T^8)$           & $\pi_B \in -4\mathbb{Q}[\pi_B]^4$ \\
					$11$ & $V_{11}(B)$                    & $B^{10}$    & 1   & 1        & $\chi_{V_{11}(B),q}(T)$               &                          \\
					$11$ & $\mid V_{11}(B)$                  & $B^5$       & 1   & 2        & $\mid \chi_{V_{11}(B),q}(T)$             &                          \\
					$11$ & $\mid V_{11}(B)$                  & $B^2$       & 1   & 5        & $\mid \chi_{V_{11}(B),q}(T)$             &                          \\
					$12$ & $\mid W_{4,1}(V_3(B/\mathbb{F}_{q^4}))$    & $B^2$       & 4   & 5        & $\mid \chi_{V_3(B/\mathbb{F}_{q^4}),{q^4}}(T^4)$  &                          \\
					$22$ & $\mid W_{2,1}(V_{11}(B/\mathbb{F}_{q^2}))$ & $B^2$       & 2   & 5        & $\mid \chi_{V_{11}(B/\mathbb{F}_{q^2}),q^2}(T^2)$ &                          \\ \bottomrule
				\end{tabular}
			\end{table}

			\section{Examples}
			\label{sec:examples}
			
			To find examples of ordinary simple geometrically split abelian varieties for each dimension considered, we used the following methods.
			\begin{enumerate}
				\item \textit{Search in LMFDB~\cite{LMFDB}}. This database contains already the information about geometrical simplicity of abelian varieties. Moreover, we can quickly find geometrically split abelian varieties with minimal splitting degree~$d$ by looking at the field of definition of endomorphisms (this information is also presented in the database). For an ordinary simple geometrically split abelian variety~$A$ the endomorphism algebra~$\operatorname{End}^0{A}$ is a field, while for~$A$ over~$\mathbb{F}_{q^d}$ the endomorphism algebra contains zero-divisors due to splitting. Therefore, for geometrically split abelian varieties there are endomorphisms defined only over~$\mathbb{F}_{q^d}$ and undefined over~$\mathbb{F}_{q}$. Thus, if we are looking for examples of geometrically split abelian variety with minimal splitting degree~$d$, we can look only for ones with full endomorphism ring defined over the field of degree~$\geq d$. However, the LMFDB contains data only for small dimensions and fields. For other cases we have to use different methods.
				
				\item
				\textit{Enumeration of Weil polynomials}. SageMath computer algebra already contains the method \texttt{weil\_polynomials} that builds a list of all possible Weil polynomials of given degree. To find examples, we enumerated Weil polynomials for small field sizes~$q$ and selected irreducible ordinary polynomials~$\chi_{A,q}$ such that~$\chi_{A,q^d}$ is reducible for some small~$d$. The polynomial~$\chi_{A,q^d}$ can be computed from~$\chi_{A,q}$ by known recurrent formulae. Factorization of polynomial over~$\mathbb{Z}$ is polynomial time procedure. An implementation of this method can be found in our source code repository.\footnote{The script is available in \texttt{examples/weil\_polynomials\_enum.sage} at \url{https://github.com/novoselov-sa/avGsplit/}}
				
				\item 
				\label{examples:method:sampling}
				\textit{Enumeration or random sampling of~$B$}. For $\dim B = 1$ we enumerated all possible elliptic curves up to isomorphism by enumerating $j$-invariants and all twists of curves for each $j$-invariant. For $\dim B \geq 2$ we sampled random hyperelliptic curves. For each randomly sampled ordinary curve we computed the characteristic polynomial~$\chi_W$ for abelian variety $W_{u,1}(V_{v}(B))$ (as required for example by Table~\ref{table:av_gsplit_dim4}) and took factors of $\chi_W$ as candidates for~$\chi_{A,q}$. Among candidates we took only ones with required minimal splitting degree.
				An implementation of this method can also be found in our source code repository.\footnote{The script is available in \texttt{examples/jacobians\_enum.sage} at \url{https://github.com/novoselov-sa/avGsplit/}}
			\end{enumerate}
			
			\subsection{Dimension $4$}
			For each case in Theorem~\ref{th:g4:frobs}, there is an explicit example in dimension~$4$. We provide them in Table~\ref{tab:g4:examples}. All of them are presented in the LMFDB.
			\begin{table}[]
				\scriptsize
				\caption{Examples of simple ordinary geometrically split dimension $4$ abelian varieties.}
				\label{tab:g4:examples}
				\begin{tabular}{lllll}
					\hline
					$d$ &
					$e$ &
					$\chi_{A,q}$ &
					$\chi_{B,q^e}$ &
					Ref. \\ \hline
					$2$ &
					$2$ &
					$T^{8} - T^{6} - T^{4} - 4 T^{2} + 2^4$ &
					$T^{4} - T^{3} - T^{2} - 4 T + 2^4$ &
					\cite[\href{https://www.lmfdb.org/Variety/Abelian/Fq/4/2/a_ab_a_ab}{4.2.a\_ab\_a\_ab}]{LMFDB} \\
					$3$ &
					$1$ &
					$T^{8} - 2 T^{7} + T^{6} + 2 T^{5} - 3 T^{4} + $ &
					$T^{4} + 2 T^{3} + 3 T^{2} + 4 T + 4$ &
					\cite[\href{https://www.lmfdb.org/Variety/Abelian/Fq/4/2/ac_b_c_ad}{4.2.ac\_b\_c\_ad}]{LMFDB} \\
					&
					&
					$+4 T^{3} + 4 T^{2} - 16 T + 16$ &
					&
					\\
					$4$ &
					$4$ &
					$T^{8} - 2 T^{7}  + 2  T^{6} - 2 T^{5} + $ &
					$T^{4} - 2 T^{3} - 15 T^{2} - 32 T + 256$ &
					\cite[\href{https://www.lmfdb.org/Variety/Abelian/Fq/4/2/ac_c_ac_b}{4.2.ac\_c\_ac\_b}]{LMFDB} \\
					&
					&
					$+T^{4} - 4 T^{3} + 8 T^{2} - 16 T + 16$ &
					&
					\\
					$5$ &
					$1$ &
					$T^{8} - T^{7} -  T^{6} + 3 T^{5} - T^{4} +$ &
					$T^{2} + T + 2$ &
					\cite[\href{https://www.lmfdb.org/Variety/Abelian/Fq/4/2/ab_ab_d_ab}{4.2.ab\_ab\_d\_ab}]{LMFDB} \\
					&
					&
					$+6 T^{3} - 4 T^{2} - 8 T + 16$ &
					&
					\\
					$5$ &
					$1$ &
					$T^{8} - 2 T^{7}  + 2 T^{6} - 5 T^{5} + $ &
					$T^{4} - T^{3} + 3 T^{2} - 2 T + 4$ &
					\cite[\href{https://www.lmfdb.org/Variety/Abelian/Fq/4/2/ac_c_af_l}{4.2.ac\_c\_af\_l}]{LMFDB} \\
					&
					&
					$+ 11 T^{4} - 10 T^{3} + 8 T^{2} - 16 T + 16$ &
					&
					\\
					$6$ &
					$2$ &
					$T^{8} - 6 T^{7} + 20 T^{6} - 48  T^{5} + 91 T^{4} -$ &
					$T^{4} - 52 T^{3}  + 2026 T^{2} - 52 \cdot 3^6 T + 3^{12}$ &
					\cite[\href{https://www.lmfdb.org/Variety/Abelian/Fq/4/3/ag_u_abw_dn}{4.3.ag\_u\_abw\_dn}]{LMFDB} \\
					&
					&
					$ - 144 T^{3} + 180 T^{2} - 162 T + 81$ &
					&
					\\
					$8$ &
					$8$ &
					$T^{8} - 4  T^{7}  + 4 T^{6} + 12 T^{5} - 40 T^{4} + $ &
					$ T^4 + 4 T^3 - 5306 T^2 + 4 \cdot 3^8 T + 3^{16}$
					&
					\cite[\href{https://www.lmfdb.org/Variety/Abelian/Fq/4/3/ae_e_m_abo}{4.3.ae\_e\_m\_abo}]{LMFDB} \\
					&
					&
					$36 T^{3} + 36 T^{2} - 108 T + 81$ &
					&
					\\
					$10$ &
					$2$ &
					$T^{8} - 5  T^{7} + 12  T^{6} - 20  T^{5} + $ &
					$T^{4} - 3 T^{3}  + 9 T^{2} - 12 T + 16$ &
					\cite[\href{https://www.lmfdb.org/Variety/Abelian/Fq/4/2/af_m_au_bd}{4.2.af\_m\_au\_bd}]{LMFDB} \\
					&
					&
					$ + 29 T^{4} - 40 T^{3} + 48 T^{2} - 40 T + 16$ &
					&
					\\
					$12$ &
					$4$ &
					$T^{8} - 4 T^{7} + 5 T^{6} + 2 T^{5} - 11 T^{4} + $ &
					$T^4 - 2 T^3 - 15 T^2 - 32 T + 256$
					&
					\cite[\href{https://www.lmfdb.org/Variety/Abelian/Fq/4/2/ae_f_c_al}{4.2.ae\_f\_c\_al}]{LMFDB} \\
					&
					&
					$4 T^{3} + 20 T^{2} - 32 T + 16$ &
					&
					\\ \hline
				\end{tabular}
			\end{table}
			However, for the case $d = 8$ the LMFDB does not contain information about~$B$. Therefore, we provide an additional example for this case.
			
			\begin{example}[$d=8$] Let~$C: y^2 = x^9 + a x^5 + b x$ be a genus~$4$ hyperelliptic curve. When the Jacobian~$\operatorname{Jac}_C$ is simple, it decomposes~\cite[\S7]{Novoselov2020} over~$\mathbb{F}_{q}[\sqrt[8]{b}]$ as~$\operatorname{Jac}_C \sim \operatorname{Jac}_{X}^2$, where~$X: y^2 = (x + 2\sqrt[8]{b})(x^4 - 4 x^2\sqrt[4]{b} + 2\sqrt{b} + a)$. Now, we can find examples of geometrically split dimension~$4$ abelian varieties with~$d=8$ by enumerating~$b \in \mathbb{F}_{q} \setminus (\mathbb{F}_{q})^8$. Recall that~$d$ is minimal, i.e. we should exclude cases where~$\operatorname{Jac}_C$ is not simple over any proper subfield of~$\mathbb{F}_{q^8}$. One such example is the following. Take~$q=17$ and~$C: y^2 = x^9 + x^5 + 3 x$. Then
				\begin{itemize}
					\item $\chi_{\operatorname{Jac}_C,q}(T) = T^8 - 8 T^7 + 24 T^6 - 32 T^5 + 32 T^4 - 32 \cdot 17 \cdot T^3 + 24 \cdot 17^2 \cdot T^2 - 8 \cdot 17^3 \cdot T + 17^4$
					\item $\chi_{\operatorname{Jac}_X,q^8}(T) = T^4 - 127228 T^3 + 12148963078 T^2 - 127228 \cdot 17^8 \cdot T + 17^{16}$
				\end{itemize}
				Polynomials $\chi_{\operatorname{Jac}_C,q}, \chi_{\operatorname{Jac}_C,q^2}$, and $\chi_{\operatorname{Jac}_C,q^4}$ are irreducible. Therefore~$\operatorname{Jac}_C$ is simple over all proper subfields of~$\mathbb{F}_{q^8}$, while over~$\mathbb{F}_{q^8}$ we have~$\chi_{\operatorname{Jac}_C,q^8} = \chi_{\operatorname{Jac}_X,q^8}^2$ as expected.
			\end{example}

			\subsection{Dimension $5$}
			By Theorem~\ref{th:av:gsplit:prime:dim}, there are two cases: $A \sim W_{5,1}(B)$ and~$A \mid V_{11}(B)$. In both cases the abelian variety~$B$ is an elliptic curve and examples can be found efficiently by enumerating $j$-invariants of elliptic curves by Method~\ref{examples:method:sampling}.
			
			\begin{example}[$A \sim W_{5,1}(B)$]
				Let $B$ be an elliptic curve defined by the equation $y^2 + x y = x^3 + (\alpha^4+\alpha)$ over a finite field~$\mathbb{F}_{2^5} = \mathbb{F}_2[\alpha] / \left< \alpha^5 + \alpha^2 + 1\right>$. Then
				\begin{enumerate}
					\item $\chi_{B,q^5}(T) = T^2 + 7 T + 32$,
					\item $\chi_{A,q} = T^{10} + 7 T^5 + 32$,
					\item $\chi_{A,q^5} = (T^2 + 7 T + 32)^5$.
				\end{enumerate}
				Therefore, $A \sim W_{5,1}(B)$ as expected.
			\end{example}
			
			\begin{example}[$A \mid V_{11}(B)$]
				Let $B$ be an elliptic curve defined by $y^2 = x^3 + 2 x^2 + 1$ over a finite field $\mathbb{F}_3$.
				\begin{enumerate}
					\item $\chi_{B,q} = T^2 + T + 3$,
					\item $\chi_{V_{11}(B),q} = (T^{10} - 6 T^9 + 14 T^8 - 7 T^7 - 46 T^6 + 133 T^5 - 138 T^4 - 63 T^3 + 378 T^2 - 486 T + 243) (T^{10} + 5 T^9 + 14 T^8 + 26 T^7 + 42 T^6 + 67 T^5 + 126 T^4 + 234 T^3 + 378 T^2 + 405 T + 243)$,
					\item $\chi_{A,q} = T^{10} - 6 T^9 + 14 T^8 - 7 T^7 - 46 T^6 + 133 T^5 - 138 T^4 - 63 T^3 + 378 T^2 - 486 T + 243$,
					\item $\chi_{A,q^{11}} = (T^2 + 67 T + 177147)^5$.
				\end{enumerate}
				Therefore, $A \mid V_{11}(B)$ as expected.
			\end{example}
			
			\subsection{Dimension~$6$}
			
			\begin{table}[H]
				\caption{Examples of simple ordinary geometrically split abelian varieties of dimension~$6$}
				\label{table:dim6:examples}
				\begin{tabular}{@{}cccccc@{}}
					\toprule
					$A$                                     & $d$  & $A/\mathbb{F}_{q^d}$ & $e$ & $\dim{B}$ & Example                                                                                                            \\ \midrule
					$W_{2,1}(B)$           & $2$  & $B^2$                & 2   & 3         & \cite[\href{https://www.lmfdb.org/Variety/Abelian/Fq/6/2/a_ah_a_ba_a_acl}{6.2.a\_ah\_a\_ba\_a\_acl}]{LMFDB}     \\
					$W_{3,1}(B)$           & $3$  & $B^3$                & 3   & 2         & \cite[\href{https://www.lmfdb.org/Variety/Abelian/Fq/6/2/a_a_e_a_a_p}{6.2.a\_a\_e\_a\_a\_p}]{LMFDB}             \\
					$V_3(B)$                                & $3$  & $B^2$                & 1   & 3         & \cite[\href{https://www.lmfdb.org/Variety/Abelian/Fq/6/2/ae_h_ag_d_ah_p}{6.2.ae\_h\_ag\_d\_ah\_p}]{LMFDB}       \\
					$\mid W_{4,1}(B)$      & $4$  & $B^2$                & 4   & 3         & Ex.~\ref{ex:dim6:d4:s2}                                                            \\
					$\mid W_{2,1}(V_3(B/\mathbb{F}_{q^2}))$ & $6$  & $B^2$                & 2   & 3         & \cite[\href{https://www.lmfdb.org/Variety/Abelian/Fq/6/2/ad_g_aj_m_ap_r}{6.2.ad\_g\_aj\_m\_ap\_r}]{LMFDB}       \\
					$V_7(B)$                                & $7$  & $B^6$                & 1   & 1         & Ex.~\ref{ex:dim6:d7:s6}                                                                                            \\
					$\mid V_7(B)$                           & $7$  & $B^3$                & 1   & 2         & Ex.~\ref{ex:dim6:d7:s3}                                                                                            \\
					$\mid V_7(B)$                           & $7$  & $B^2$                & 1   & 3         & \cite[\href{https://www.lmfdb.org/Variety/Abelian/Fq/6/2/f_j_h_d_b_ab}{6.2.f\_j\_h\_d\_b\_ab}]{LMFDB}           \\
					$\mid W_{8,1}(B)$      & $8$  & $B^2$                & 8   & 3         & $-$                                                                                                                \\
					$\mid V_9(B)$                           & $9$  & $B^2$                & 1   & 3         & \cite[\href{https://www.lmfdb.org/Variety/Abelian/Fq/6/2/ag_p_at_g_bb_acj}{6.2.ag\_p\_at\_g\_bb\_acj}]{LMFDB}
					
					\\
					$\mid W_{3,1}(V_3(B/\mathbb{F}_{q^3}))$ & $9$  & $B^2$                & 3   & 3         & Ex.~\ref{ex:dim6:d9:e3} \\
					$\mid W_{4,1}(V_3(B/\mathbb{F}_{q^4}))$ & $12$ & $B^2$                & 4   & 3         & $-$                                                                                                                \\
					$\mid V_{13}(B)$                        & $13$ & $B^3$                & 1   & 2         & Ex.~\ref{ex:dim6:d13}                                                                                              \\
					$\mid W_{2,1}(V_7(B/\mathbb{F}_{q^2}))$ & $14$ & $B^2$                & 2   & 3         & $-$                                                                                                                \\
					$\mid W_{2,1}(V_9(B/\mathbb{F}_{q^2}))$ & $18$ & $B^2$                & 2   & 3         & \cite[\href{https://www.lmfdb.org/Variety/Abelian/Fq/6/2/ag_p_av_y_abn_cn}{6.2.ag\_p\_av\_y\_abn\_cn}]{LMFDB} \\
					$\mid W_{6,1}(V_3(B/\mathbb{F}_{q^6}))$ & $18$ & $B^2$                & 6   & 3         & 
					$-$
					\\
					$\mid W_{3,1}(V_7(B/\mathbb{F}_{q^3}))$ & $21$ & $B^3$                & 3   & 2         & $-$                                                                                                                \\ \bottomrule
				\end{tabular}
			\end{table}
			
			\begin{example}[$d = 4$, $\dim{B} = 3$]
				\label{ex:dim6:d4:s2} To find an example for $d = 4$, we enumerated all irreducible Weil polynomials~$f(T)$ of degree $6$ over the field~$\mathbb{F}_{3^4}$ and checked that $f$ is irreducible, ordinary, and in the number field~$K = \mathbb{Q}(\pi_B) = \mathbb{Q}[T] / \left< f(T) \right>$ the element~$-T/4$ is a fourth power. If this condition is satisfied, we have $\pi_B = -4 \mathbb{Q}[\pi_B]^4$ as required.
				For the enumeration, we used the SageMath function~$\texttt{weil\_polynomials}$, which constructs all Weil polynomials of given degree. There are $12104353$ Weil polynomials over~$\mathbb{F}_{3^4}$. Only a few of them satisfy the required conditions. For example, take
				\[
				f = T^6 - 44 T^5 + 872 T^4 - 10022 T^3 + 70632 T^2 - 288684 T + 531441.
				\]
				Since Weil polynomials of degree~$6$ define abelian varieties of dimension~$3$, we can take $\chi_{B,q}(T) = f$.
				Then
				$\chi_{W_{4,1}(B),q}(T) = \chi_{B,q^4}(T^4) = (T^{12} - 2 T^{11}+2 T^{10}-8 T^9-8 T^8+38 T^7-28 T^6+114 T^5-72 T^4-216 T^3+162 T^2-486 T+729) (T^{12}+2 T^{11}+2 T^{10}+8 T^9-8 T^8-38 T^7-28 T^6-114 T^5-72 T^4+216 T^3+162 T^2+486 T+729)$. By choosing the first one from these factors as $\chi_{A,q}(T)$, we obtain an example of a geometrically split abelian variety. The polynomial~$\chi_{A,q^2}(T) = T^{12} - 44 T^{10} + 872 T^8 - 10022 T^6 + 70632 T^4 - 288684 T^2 + 531441$ is irreducible and~$\chi_{A,q^2}(T) = \chi_{B,q}(T)^2$. Therefore, the minimal splitting degree is~$4$ as required.
			\end{example}
			
			\begin{example}
				\label{ex:dim6:d7:s6} Take an elliptic curve~$E: y^2 = x^3 + x^2 + 2$ over~$\mathbb{F}_{3}$. Then $V_7(E)$ has dimension~$6$, its characteristic polynomial~$\chi_{V_7(E),q}(T) = T^{12} + T^{11} - 2 T^{10} - 5 T^9 + T^8 + 16 T^7 + 13 T^6 + 48 T^5 + 9 T^4 - 135 T^3 - 162 T^2 + 243 T + 729$ is irreducible and~$\chi_{V_7(E),{q^7}}(T) = (T^2 + 83 T + 2187)^6$. Therefore, $A = V_7(E)$ is a geometrically split simple ordinary abelian variety of dimension~$6$.
			\end{example}
			
			\begin{example}
				\label{ex:dim6:d7:s3} Take a hyperelliptic curve~$C: y^2 = x^5 + 110 x^4 + 89 x^3 + 180 x^2 + 39 x + 120$ over~$\mathbb{F}_{193}$. The characteristic polynomial of~$\operatorname{Jac}_C$ is $T^4 - 62 T^2 + 37249$.
				Then $V = V_7(\operatorname{Jac}_C)$ has dimension~$12$ and its characteristic polynomial~$\chi_{V,q}(T)$ factors into two irreducible polynomials:
				\begin{itemize}
					\item $T^{12} - 473 T^{10} + 78577 T^8 - 8619353 T^6 + 78577 \cdot 193^2 T^4 - 473 \cdot 193^4 T^2 + 193^6$,
					\item $T^{12} + 535 T^{10} + 141073 T^8 + 28927639 T^6 + 141073 \cdot 193^2 T^4 + 535 \cdot 193^4 T^2 + 193^6$.
				\end{itemize}
				Take as~$\chi_{A,q}$ the first polynomial. Then \[
				\chi_{A,q^7}(T) = (T^4 + 18036094842549442 T^2 + 193^{14})^3.\]
				Therefore, $A$ is an example of a geometrically split abelian variety of dimension~$6$.
			\end{example}
			
			\begin{example}[$\dim A = 6$, $d = 9$, $B$ is primitive over~$\mathbb{F}_{q^3}$]
				\label{ex:dim6:d9:e3}
				This example was found by enumerating random hyperelliptic curves over~$\mathbb{F}_{3^3}$. For each curve, we computed the characteristic polynomial of its Jacobian~$\operatorname{Jac}_{C}$. We then used this polynomial to compute the characteristic polynomial of the abelian variety $W_{3,1}(V_3(B/\mathbb{F}_{3^3}))$ and took the irreducible factors from its factorization. For each irreducible factor, which we considered a candidate for~$\chi_{A,q}$, we checked that $\chi_{A,q^3}$ is irreducible, so $A$ remains simple over~$\mathbb{F}_{q^3}$. After examining approximately~$1{,}000$ curves, we found the following example.
				\begin{itemize}
					\item $B = \operatorname{Jac}_C$ for $C$ defined over~$\mathbb{F}_{3^3}$ by the equation $y^2 = x^8 + (2 \alpha^2 + 2 \alpha) x^7 + 2 x^6 + (\alpha^2 + 2 \alpha + 2) x^5 + (2 \alpha^2 + \alpha + 1) x^4 + (\alpha^2 + 2 \alpha + 1) x^3 + (2 \alpha^2 + 2 \alpha + 2) x^2 + (2 \alpha^2 + 2) x + 2 \alpha^2 + 1$
					\item $\mathbb{F}_{3^3} = \mathbb{F}_3(\alpha)/\left<\alpha^3 + 2\alpha + 1\right>$
					\item $\chi_{B,q^3} = T^6 - 3 T^5 + 3 T^4 + 161 T^3 + 81 T^2 - 2187 T + 19683$
					\item $\chi_{A,q}(T) = T^{12} - 9 T^{11} + 45 T^{10} - 161 T^9 + 450 T^8 - 1026 T^7 + 1945 T^6 - 3078 T^5 + 4050 T^4 - 4347 T^3 + 3645 T^2 - 2187 T + 729$
				\end{itemize} 
			\end{example}
			
			\begin{example}[$\dim A = 6$, $d = 13$]
				\label{ex:dim6:d13} Take a hyperelliptic curve~$C: y^2 = x^5 + x^3 + x + 2$ over~$\mathbb{F}_{3}$. The characteristic polynomial of~$\operatorname{Jac}_C$ is $T^4 - 3 T^3 + 5 T^2 - 9 T + 9$.
				Then $V = V_{13}(\operatorname{Jac}_C)$ has dimension~$12$ and its characteristic polynomial~$\chi_{V,q}(T)$ factors into~$4$ irreducible polynomials. One of these is~$\chi_{A,{q}} = T^{12} - 9 T^{11} + 42 T^{10} - 131 T^9 + 308 T^8 - 601 T^7 + 1067 T^6 - 1803 T^5 + 2772 T^4 - 3537 T^3 + 3402 T^2 - 2187 T + 729$. Then $
				\chi_{A,q^{13}}(T) = (T^4 + 3975 T^3 + 7045979 T^2 + 6337433925 T + 2541865828329)^3$.
				Therefore, $A$ is an example of a geometrically split abelian variety of dimension~$6$.
			\end{example}

		\end{document}